\documentclass[a4paper]{article}
\usepackage[all]{xy}
\usepackage[latin1]{inputenc}        
\usepackage{graphics,graphicx}
\usepackage{amsfonts,amssymb,amsmath,color,mathrsfs,amstext,bm}
\usepackage{amsbsy, amsopn, amscd, amsxtra, amsthm,authblk}
\usepackage{enumerate,algorithmic,algorithm,adjustbox}
\usepackage{upref}
\usepackage{geometry}
\usepackage[mathlines]{lineno}
\usepackage{float}
\usepackage{yhmath}
\usepackage{booktabs}
\usepackage{subcaption}
\usepackage{multirow}
\usepackage{makecell}
\usepackage{tikz}
\usetikzlibrary{arrows.meta,backgrounds,calc,fit}
\usepackage[colorlinks,
            linkcolor=red,
            anchorcolor=red,
            citecolor=red
            ]{hyperref}
\usepackage[nameinlink,noabbrev]{cleveref}

\numberwithin{equation}{section}
\numberwithin{figure}{section}
\numberwithin{table}{section}

\def\R{\mathbb{R}}
\def\P{\mathbb{P}}
\def\E{\mathbb{E}}
\def\cL{\mathcal{L}}
\def\cD{\mathcal{D}}
\def\cI{\mathcal{I}}
\def\cP{\mathcal{P}}

\def\cQ{\mathcal{Q}}

\def\cH{\mathcal{H}}

\def\cF{\mathcal{F}}
\def\cB{\mathcal{B}}
\def\cK{\mathcal{K}}

\def\tr{\mathrm{tr}}
\def\l{\langle}
\def\r{\rangle}

\newcommand{\dd}{\,\mathrm d}
\newcommand{\Id}{I_d}
\newcommand{\one}{\mathbf 1}

\newcommand{\Law}{\operatorname{Law}}
\newcommand{\Wtwo}{W_2}

\newtheorem{theorem}{Theorem}[section]
\newtheorem{lemma}{Lemma}[section]
\newtheorem{proposition}{Proposition}[section]

\newtheorem{corollary}{Corollary}[section]

\newtheorem{assumption}{Assumption}[section]
\crefname{assumption}{assumption}{assumptions}
\Crefname{assumption}{Assumption}{Assumptions}

\makeatletter
\@mparswitchfalse
\makeatother

\title{Quantitative Convergence of the Spherical Brownian Motion collisional particle system for the Maxwellian Landau Equation}
\author[a]{Jos\'e Antonio Carrillo\thanks{
    E-mail: \texttt{jose.carrillo@maths.ox.ac.uk}.}}
\author[b]{Ruizhi Huang\thanks{
    E-mail: \texttt{ruizhi.huang36@gmail.com}.}}
\author[b,c]{Lei Li\thanks{ E-mail: \texttt{leili2010@sjtu.edu.cn}.}}

\affil[a]{Mathematical Institute, University of Oxford,
    Woodstock Road, Oxford OX2 6GG, United Kingdom.}
\affil[b]{School of Mathematical Sciences,
    Shanghai Jiao Tong University, Shanghai 200240, China.}
\affil[c]{Institute of Natural Sciences, MOE-LSC, Shanghai Jiao Tong University, Shanghai 200240, China.}

\date{}

\begin{document}

\maketitle

\begin{abstract}
    In this paper, we establish the quantitative convergence of the spherical Brownian motion (SBM) collisional particle system for the Maxwellian Landau equation via Wasserstein distances. When used as numerical method, it enjoys linear complexity $O(N)$ in particle number $N$, and conserves the momentum and energy almost surely. Moreover, only random grouping and sampling are needed for simulations without further numerical discretization.  In the first part, we identify the large particle mean field limit of the system, and the Wasserstein convergence to the limit is justified through a coupling similar to the recent work for the Nanbu type Monte Carlo method \cite{BorghiPareschi2025}. The key difference is that the particles in our coupled shadow system are not i.i.d. due to the perfect matching. In the second part, the convergence of the mean field dynamics to the Landau solution is justified by the relative entropy with a bound $O(\tau^2)$, being $\tau$ the stepsize. This results in an $O(\tau)$ bound for the Wasserstein distance error. A crucial technical tool in the convergence in relative entropy is to provide suitable bounds on the derivatives of the logarithmic density or score function in the mean field dynamics. The combination of these estimates leads to the full convergence of the squared Wasserstein distance with order $O(N^{-1/2}+\tau^2)$ for the numerical scheme in $d=3$.
\end{abstract}

\section{Introduction}

The Landau equation is a fundamental kinetic model for collisional relaxation in plasmas, where charged particles interact through Coulomb forces and grazing collisions predominate \cite{Landau1936,Villani2002,CarrilloDelgadinoWu2022Grazing}. Its collision operator also enters spatially inhomogeneous Vlasov--Landau models, which couple collisions with particle transport and self-consistent fields \cite{DimarcoLiPareschiYan2015,BailoCarrilloHu2024}. In the spatially homogeneous setting considered here, the equation reads as follows:
\begin{equation}\label{eq:landau}
    \partial_t f=\cQ(f,f)
    :=\nabla_v\cdot\int_{\R^d} A_\gamma(v-w)
    \bigl(f(t,w)\nabla_v f(t,v)-f(t,v)\nabla_w f(t,w)\bigr)\dd w.
\end{equation}
Here, \(d\ge2\) is the dimension of the velocity, \(\Lambda>0\) is a fixed constant, and \(f(t,\cdot)\) is the probability density on \(\R^d\).  The interaction kernel matrix is given by
\begin{equation}\label{eq:kernel}
    A_\gamma(z)=\Lambda|z|^\gamma\bigl(|z|^2I_d-z\otimes z\bigr),\qquad -3\le\gamma\le1.
\end{equation}
In this work, we focus on Maxwell molecules, corresponding to \(\gamma=0\). The conservation of mass, momentum, and kinetic energy, together with entropy dissipation \cite{CarrilloDelgadinoDesvillettesWu2024} are the important physical properties featured by this model.
They also provide important criteria for assessing numerical methods for approximating its solutions.

Typical numerical methods include fast spectral schemes \cite{PareschiRussoToscani2000} and finite element methods \cite{HirvijokiAdams2017,AdamsEtAl2017}. These approaches often face the difficulty introduced by the nonlocality, nonlinearity and degeneracy of the models.  A more convenient way is to seek the particle approximation of the Landau models.
Typical stochastic particle approximations to Landau and Landau-type equations include conservative Kac-type diffusion systems \cite{Carrapatoso2016} and interacting diffusion systems based on the McKean--Vlasov formulation \cite{Fournier2009}. These belong to the probabilistic models and have been studied recently by several groups, using coupling techniques \cite{FGHard}, relative entropy methods \cite{CarrilloGuoJabin2024,CFGJW}, and some hierarchy type techniques \cite{FengWang2026Kac}.
These particle systems, however, are not very suitable for numerical simulations because they usually cost $O(N^2)$ per time step if applied as numerical methods. Some existing numerical particle methods include the Monte Carlo collision algorithms \cite{BobylevNanbu2000,MedagliaPareschiZanella2024}, and deterministic particle methods based on a regularized entropy \cite{CarrilloHuWangWu2020,CarrilloHuVanFleet2024,BailoCarrilloHu2024,BailoCarrilloMedagliaZanella2025} and the random batch version \cite{CarrilloJinTang2022}.  The particle methods based on regularized entropy are numerical methods aiming to approximate the densities and do not keep the particle microscopic collision properties. The Monte Carlo methods in \cite{BobylevNanbu2000,MedagliaPareschiZanella2024} have $O(N)$ complexity but their approximations are based on the approximation to the Boltzmann equations instead of on the Landau operator directly.
The grazing regime is intrinsically the diffusion regime of elastic collision, and some particle systems have been proposed for the Landau and Landau type equations \cite{du2024collision,du2025structure} building on this observation. Their advantage is that these methods are direct approximations to Landau collision operator instead of approximating the Boltzmann operator in grazing approximation.
We remark that the particle models in \cite{du2024collision,du2025structure} look similar to the random batch type methods  \cite{JinLiLiu2020,JinLiLiu2021Disparate}, however they are derived directly from the pairwise collision in the grazing regime, and thus, they are not approximations to the fully-interacting collisional $N$ particle system.

The diffusion approximation, as crucially observed in \cite{du2025structure}, makes the relative velocity change  due to grazing collisions similar to a spherical Brownian motion. In \cite{du2025structure},  the Spherical Brownian Motion collisional particle system (SBM) is then proposed for the approximation of the Landau equation, which organizes the particles into a uniform random perfect matching at each time step. Within each matched pair, the center-of-mass velocity and the magnitude of the relative velocity remain fixed, while the direction of the relative velocity evolves by spherical Brownian motion. \Cref{alg:sbm} recalls the SBM scheme introduced in \cite[Algorithm~2]{du2025structure}, written with the normalization of \eqref{eq:landau}\footnote{Equation~\eqref{eq:landau} omits the factor \(1/2\) in \cite[Eq.~(1)]{du2025structure}. For the same \(\Lambda\), this doubles the collision generator and changes the SBM time-scaling coefficient from \(4\Lambda|z|^\gamma\) to \(8\Lambda|z|^\gamma\). If the two pre-collision velocities coincide, you leave them unchanged and skip the remaining collision calculations.}.  When used as a numerical method, the model enjoys linear complexity $O(N)$ per time step and conserves the momentum and energy almost surely. More importantly, only random grouping and sampling from the spherical Brownian motions are needed for the simulation. We emphasize that the model is not only a numerical method, but also plays as physical approximation to the Landau collision. The SBM model enjoys a lot of similarity with the direct simulation Monte Carlo (DSMC) methods for Boltzmann equations. Classical DSMC methods include the one-sided Nanbu scheme \cite{Nanbu1980} and the Bird scheme \cite{Bird1994}. In the Nanbu-type,
each particle independently samples a partner and collision parameters, and selecting a partner does not require a reciprocal update. For such methods for elastic collisions, momentum and energy are conserved in expectation, but not generally pathwise. The convergence in Wasserstein distance of the Nanbu-type of DSMC has recently been proved by Borghi and Pareschi \cite{BorghiPareschi2025}. The SBM model is more similar to the Bird-type, particles are updated in matched pairs using shared randomness.

\begin{algorithm}[H]
    \caption{Spherical Brownian motion (SBM) scheme.}
    \label{alg:sbm}
    \begin{algorithmic}[1]
        \REQUIRE Particle number $N$ (even), initial velocity of particles
        $V_0^i$ $(i=1,\ldots,N)$, time step $\tau$, terminal time $T=K\tau$.
        \FOR{$k=0:K-1$}

        \STATE Randomly divide $N$ particles into $N/2$ pairs, and let $\theta_k$ denote the matching function.

        \STATE For each pair, compute
        \[
            z=V_k^i-V^{\theta_k(i)}_k,
            \qquad
            s=V^i_k+V^{\theta_k(i)}_k,
            \qquad
            e_z=\frac{z}{\lvert z\rvert}.
        \]

        \STATE  Calculate the time-scaling coefficient
        \[
            \lambda=8\Lambda\lvert z\rvert^\gamma.
        \]

        \STATE Simulate a standard SBM on $\mathbb{S}^{d-1}$
        starting from $e_z$ with a sampling time $\lambda\tau$.
        Denote the result as $e_z'$.

        \STATE Compute the post-collision velocity difference
        $z'=\lvert z\rvert e_z'$, and update the velocities
        \[
            V^i_{k+1}=\frac{s+z'}{2},
            \qquad
            V^{\theta_k(i)}_{k+1}=\frac{s-z'}{2}.
        \]
        \ENDFOR

        \STATE \textbf{Return} the velocities at terminal time:
        $V^i_{T/\tau}$ $(i=1,\ldots,N)$.
    \end{algorithmic}
\end{algorithm}

The study of the convergence of the SBM model to the Landau equation for Maxwellian molecules not only justifies the correctness of the method when used as numerical methods, but also supports the micro-interpretation of Landau model.
This paper then aims to provide the convergence of the SBM model.  The analysis involves two steps. The first is to justify the large particle limit $N\to\infty$. The second is to prove the convergence of the limit scheme to the Landau equation.

In the first step, our key observation is that the spherical Brownian time scale is independent of the relative speed for Maxwell molecules, so a common random rotation gives a linear representation of the collision update. We exploit this representation to construct a coupling in similar spirit as Borghi and Pareschi \cite{BorghiPareschi2025}. As commented above,  the resulting dependence in the updates of matched particles must be accounted for in the convergence analysis. This becomes the key difference from \cite{BorghiPareschi2025}, so we have to make extra effort to deal with this dependence.
This construction then justifies the convergence as $N\to\infty$ to a limit nonlinear dynamics. In particular, our large-particle result, \Cref{thm:main}, gives a quantitative convergence rate to the limit SBM model. Let $\tau$ be the time step and let $T=K\tau$ be the computational time. For even \(N\), independent initial velocities with common law \(f_0\), and a finite moment \(M_q(f_0):=\int |v|^q f_0(\dd v)\) of some order \(q>4\), we prove the empirical measure $f_k^N$ at time point $t_k$ converges to the mean field distribution $f_k$ under the euclidean Wasserstein distance $W_2$ with the following rate
\[
    \max_{0\le k\le K}\E W_2^2(f_k^N,f_k)
    \le C(T,d,q,\Lambda)M_q(f_0)^{2/q}\varepsilon_d(N),
\]
where
\[
    \varepsilon_d(N)=
    \begin{cases}
        N^{-1/2},           & 2\le d<4, \\
        N^{-1/2}[1+\log N], & d=4,      \\
        N^{-2/d},           & d>4.
    \end{cases}
\]
In particular, the expected squared \(W_2\) error is bounded by a multiple of \(N^{-1/2}\) in dimension three.
We establish a covariance bound for bounded shadow observables and combine it with the transport estimates of \cite{FournierGuillin2015} to obtain the same particle-number rates as for i.i.d. empirical measures.

For the second limit, we compare the nonlinear mean-field algorithm with the Landau solution through relative entropy. The analysis draws on logarithmic gradient and Hessian estimates for the reference density, as used in the Landau--Maxwellian particle analysis \cite{CFGJW}. We use these tools to compare the one-particle mean-field evolution directly with the Landau solution. This approach belongs to a broader entropy and Fisher-information framework for the Landau equation \cite{GuillenSilvestre2025,tabary2026weak}.
Let \(\widetilde f(t,\cdot)\) denote a natural continuous interpolation of the mean-field sequence \(f_k\) specified below. Under certain regularity assumptions specified below, we show in \Cref{thm:step2-entropy} that
\[
    \sup_{0\le t\le T}\cH\bigl(\widetilde f(t,\cdot)\mid f(t,\cdot)\bigr)
    \le C_T\tau^2,
\]
with the relative entropy between two probability densities $f$ and $g$ defined as
\[
    \cH(f\mid g):=\int_{\R^d}f\log\frac{f}{g}\,\dd v.
\]
Its proof combines a quantitative estimate of within-step correlations, propagation of density-derivative bounds, and relative coercivity. These two estimates yield a comparison between the particle system and the Landau solution under additional assumptions. In dimension three, they give the estimate
\[
    \max_{0\le k\le K}\E W_2^2\bigl(f_k^N,f(t_k,\cdot)\bigr)
    \le C_T\bigl(N^{-1/2}+\tau^2\bigr).
\]

The rest of the paper is organized as follows. \Cref{sec:setup} introduces the nonlinear mean-field algorithm, develops the rotation representation used in the analysis, and presents the main convergence results. \Cref{sec:large_particle} proves the time-step-uniform large-particle estimate. \Cref{sec:entropic} proves the relative-entropy estimate and the resulting error bounds for the particle scheme. The appendices provide the uniqueness argument for the two-particle forward equation, logarithmic regularity estimates for the reference solution, and the Maxwell relative-coercivity estimate.

\section{Setup, notations and main results}\label{sec:setup}

In this section, we introduce the notation and mathematical setup, establish preliminary results, and state the main convergence results.

In fact, the large particle limit of the SBM model is not hard to identify. Similar as in \cite{JinLi2022MeanFieldRBM}, one can show that as $N\to\infty$, the limiting dynamics of a focused particle (say particle $1$) is given as in \Cref{alg:sbmlimit}. The inuition is that if the particle number $N\to\infty$, the chosen particle to pair with the focused one will be independent, with the same law.
We note that, however, the proof in \cite{JinLi2022MeanFieldRBM} holds for the limit $N\to\infty$
with any fixed step number $k$ and $\tau$. It cannot be uniform for all $k$ with $k\tau\le T$ and thus the convergence as time step goes to zero is not established.

\begin{algorithm}[H]
    \caption{Limit spherical Brownian motion (SBM) scheme.}
    \label{alg:sbmlimit}
    \begin{algorithmic}[1]
        \REQUIRE Initial velocity of particles
        $V_0$ , time step $\tau$, terminal time $T=K\tau$.
        \FOR{$k=0:K-1$}

        \STATE Draw an independent copy of the current particle, denoted by $V_k^*$

        \STATE Compute
        \[
            z=V_k-V^{*}_k,
            \qquad
            s=V_k+V^{*}_k,
            \qquad
            e_z=\frac{z}{\lvert z\rvert},
            \qquad
            \lambda=8\Lambda\lvert z\rvert^\gamma.
        \]

        \STATE Simulate a standard SBM on $\mathbb{S}^{d-1}$
        starting from $e_z$ with a sampling time $\lambda\tau$.
        Denote the result as $e_z'$.

        \STATE Compute the post-collision velocity difference
        $z'=\lvert z\rvert e_z'$, and update the velocity
        \[
            V_{k+1}=\frac{s+z'}{2},
        \]
        \ENDFOR

        \STATE \textbf{Return} the velocities at terminal time:
        $V_{T/\tau}$.
    \end{algorithmic}
\end{algorithm}

We may write down the evolution of the law of $V_k$ in \Cref{alg:sbmlimit}.   Introduce the coordinates
\begin{equation}\label{eq:center_relvel}
    c:=\frac{v+v_*}{2}, \qquad \omega:=v-v_*.
\end{equation}
Since \(v=c+\omega/2\) and \(v_*=c-\omega/2\), the chain rule gives
\[
    \nabla_\omega=\frac12(\nabla_v-\nabla_{v_*}).
\]
Here, differentiation with respect to \(\omega\) is taken at fixed \(c\).
Then, $A(\omega) = \Lambda|\omega|^{\gamma+2}\Pi(\omega)$ and $\Pi(\omega) = I_d - \frac{\omega\otimes\omega}{|\omega|^2}$.
Define
\begin{gather}
    t_k=k\tau.
\end{gather}
The law $f_k:=\Law(V_k)$ is given by \Cref{alg:sbm-mean-field}. In \Cref{subsec:limitdynamics}, we will clarify why \Cref{alg:sbm-mean-field} gives the evolution of the law of $V_k$
in \Cref{alg:sbmlimit}. Intuitively, this intermediate mean field SBM particle method is intimately related to the lifting particle dynamics or the 2-particle collision in the Landau hierarchy used in \cite{GuillenSilvestre2025,carrillo2025fisherinformationdecaykac} for the analysis of the Landau equation.

\begin{algorithm}[H]
    \caption{Mean field dynamics of the SBM particle method}
    \label{alg:sbm-mean-field}
    \begin{algorithmic}
        \STATE Take $\tau>0$, $t_k=k\tau$.
        \FOR{$k=0,1,2,\ldots$}
        \STATE Set $F_k=f_k^{\otimes 2}$
        \STATE Evolve the linear equation
        \[
            \partial_t F
            =4\nabla_{\omega}\cdot
            \bigl(A(\omega)\nabla_{\omega}F\bigr)
            \mathrel{:=}\mathcal{L}^{*}F,
        \]
        to get $F_{k+1}^-=e^{\tau \mathcal{L}^*}F_k$.
        \STATE Set
        $\displaystyle
            f_{k+1}=\pi_1 F_{k+1}^{-}=\int F_{k+1}^{-}(x,y)\,dy
        $
        \ENDFOR
    \end{algorithmic}
\end{algorithm}

We denote the empirical measure of the particle system by
\begin{equation}\label{eq:fkNdef}
    f_k^N:= \frac{1}{N} \sum_{i=1}^N \delta_{V_k^i}.
\end{equation}
In the first step, we aim to show that $f_k^N$ converges to $f_k$ as $N\to \infty$ uniformly for all $k,\tau$ with $k\tau\le T$. Moreover, in the second step, we aim to show the convergence of \Cref{alg:sbm-mean-field} to the Landau equation
as $\tau\to 0$. These will then establish the convergence of the SBM model to the Landau equation for the Maxwell case $\gamma=0$.

\subsection{Rewrite the SBM scheme with Maxwell rotation}\label{subsec:sbmbasics}

We first rewrite the definition of $V_k^i$. Consider Brownian rotations $\{\mathcal{R}_s\}_{s\geq 0}$ defined by the following Stratonovich equation
\begin{equation} \label{eq:SO-Brownian-SDE}
    \dd \mathcal{R}_s = \sum_{1\leq a < b\leq d} J_{ab} \mathcal{R}_s\circ\dd B_s^{ab},\quad \mathcal{R}_0=I_d,
\end{equation}
where $J_{ab}:=e_ae_b^{\mathsf T}-e_be_a^{\mathsf T}$ with $e_i$ being standard unit vectors and $s\mapsto B_s^{ab}\in \R$ are independent standard 1D Brownian motions. To relate this rotation process to the usual projection representation of spherical Brownian motion, fix \(e\in\mathbb S^{d-1}\) and set \(n_s:=\mathcal R_s e\). The Stratonovich chain rule gives
\[
    \mathrm dn_s
    =\sum_{a<b}J_{ab}n_s\circ\mathrm dB_s^{ab},
    \qquad n_0=e.
\]
For every \(n\in\mathbb S^{d-1}\), since  $\sum_{a<b}J_{ab}^2=-(d-1)I_d$,
$ \sum_{a<b}(J_{ab}n)(J_{ab}n)^{\mathsf T}
    =I_d-nn^{\mathsf T}$, the corresponding It\^o equation for $n_s$ is
\[
    \mathrm dn_s
    =\sum_{a<b}J_{ab}n_s\,\mathrm dB_s^{ab}
    -\frac{d-1}{2}n_s\,\mathrm ds,
\]
and generates the same law as
\[
    \mathrm d\widetilde n_s
    =(I_d-\widetilde n_s\widetilde n_s^{\mathsf T})
    \,\mathrm d\tilde{B}_s
    -\frac{d-1}{2}\widetilde n_s\,\mathrm ds=(I_d-\widetilde n_s\widetilde n_s^{\mathsf T})
    \circ\mathrm d\tilde{B}_s,
\]
where \(\tilde{B}_s\) is a standard \(d\)-dimensional Brownian motion. This is the spherical Brownian motion with generator \(\frac12\Delta_{\mathbb S^{d-1}}\)\cite[Eq.~(22)]{du2025structure}, and uniqueness in law implies that \((\mathcal R_s e)_{s\geq0}\) has the same law as the solution of the projection equation\cite[Chapter~32]{Kallenberg2021}.
We use the rotation representation because it allows the same random rotation to act on different relative-velocity vectors in the coupling constructed below.

Define
\[
    P_\tau:=\mathcal{R}_{8\Lambda \tau},
\]
hence we have $e_z' = P_\tau e_z$ in \Cref{alg:sbm}.
For a rotation matrix $P$, denote
\[
    A_P:=\frac{I_d + P}{2},\qquad B_P:=\frac{I_d-P}{2}.
\]
For every matched unordered pair $\{i,j\}, j=\theta_k(i)$, at step $k$, take an independent copy $P_k^{\{i,j\}}$ of $P_\tau$ and set
\[
    P_k^i = P_k^j = P_k^{\{i,j\}}.
\]
Then we have
\begin{align}\label{eq:finite-update}
    V_{k+1}^i & = \frac{V_k^i+V_k^j}{2} + P_k^i\left(\frac{V_k^i-V_k^j}{2}\right) \\
              & =A_{P_k^i}V_k^i + B_{P_k^i}V_k^j.
\end{align}

We can calculate that
\[
    \E P_\tau = q_\tau I_d,\qquad q_\tau= e^{-4\Lambda(d-1)\tau}.
\]
Consequently, with
\begin{equation}
    a_\tau:=\frac{1+q_\tau}{2}, \qquad b_\tau:=\frac{1-q_\tau}{2},
    \label{eq:a-b}
\end{equation}
we have
\begin{equation}
    a_\tau+b_\tau=1, \qquad
    0\leq b_\tau\leq2\Lambda(d-1)\tau,
    \label{eq:b-small}
\end{equation}
and the exact matrix identities
\begin{equation}
    \E(A_{P_\tau}^{\mathsf T}A_{P_\tau})=a_\tau\Id,\quad\E(B_{P_\tau}^{\mathsf T}B_{P_\tau})=b_\tau\Id,\quad\E(A_{P_\tau}^{\mathsf T}B_{P_\tau})=0.
    \label{eq:rotation-identities}
\end{equation}
In particular, for deterministic $x,y\in\R^d$,
\begin{equation}
    \E\lvert A_{P_\tau}x+B_{P_\tau}y\rvert^2
    =a_\tau|x|^2+b_\tau|y|^2.
    \label{eq:quadratic-identity}
\end{equation}
The coefficient $b_\tau=O(\tau)$ is the mechanism that prevents the final constant from depending on $\tau$.

\subsection{The Maxwell mean field recursion}\label{subsec:limitdynamics}

For a probability measure $\mu\in\mathcal P(\R^d)$, define
\[
    M_q(\mu):=\int_{\R^d}|v|^q\mu(\dd v), \qquad q>0.
\]
We denote $\mathcal{P}_2(\R^d)$ to be the set of all probability measures with finite second moment, equipped with the Wasserstein-2 distance.

According to \Cref{alg:sbmlimit}, for probability measure $\mu$, we can define
\begin{equation}
    S_\tau\mu:=\Law(A_{P_\tau}X+B_{P_\tau}X^*),\qquad X,X^*\overset{\mathrm{i.i.d.}}{\sim}\mu,\quad P_\tau\perp(X,X^*). \label{eq:Sh}
\end{equation}
We first verify that the mean-field recursion in \Cref{alg:sbm-mean-field} for Maxwell molecules satisfies
\begin{equation}\label{eq:rho_recursion}
    f_{k+1}=S_\tau f_k.
\end{equation}
In particular, take \(\mu=f_k\), and take the rotation process \(\mathcal R\), starting from \(I_d\), independently of this initial pair. Set
\[
    C:=\frac{X+X^*}{2}, \qquad \Omega_0:=X-X^*.
\]
For the local time \(s\in[0,\tau]\), define
\begin{equation}\label{eq:candidate_pair}
    X_s:=C+\frac12\mathcal R_{8\Lambda s}\Omega_0, \qquad X_s^*:=C-\frac12\mathcal R_{8\Lambda s}\Omega_0, \qquad \nu_s:=\Law(X_s,X_s^*).
\end{equation}
We now verify that \((\nu_s)_{0\leq s\leq \tau}\) solves the initial-value problem with $s=t-t_k$
\begin{equation}\label{eq:pair_initial_value_problem}
    \begin{cases}
        \partial_t F_t=\mathcal L^*F_t :=4\nabla_\omega\cdot\bigl(A(\omega)\nabla_\omega F_t\bigr), & t\in(t_k,t_{k+1}), \\
        F_{t_k}=f_k^{\otimes2}.
    \end{cases}
\end{equation}
We denote its terminal value by
\[
    F_{k+1}^-:=\lim_{s\uparrow \tau}F_{t_k+s}\quad\text{in } \mathcal{P}_2(\mathbb R^{2d}).
\]

\begin{lemma}\label{lem:pair-law}
    Suppose that $f_k \in \mathcal{P}_2$.  Then, the law $\nu_s=\Law(X_s, X_s^*)$ is the unique solution in $C([0, \tau]; \mathcal{P}_2(\R^{2d}))$ to \Cref{eq:pair_initial_value_problem} in the weak sense for probability measures, and consequently \Cref{eq:rho_recursion}
    holds. Moreover,
    \[
        M_2(f_{k+1})=M_2(f_k)
    \]
\end{lemma}

\begin{proof}
    We now verify that the joint law \(\nu_s\) of \(Y_s:=(X_s,X_s^*)\) satisfies the prescribed forward equation. We first derive its It\^o equation from the rotation representation in \Cref{eq:candidate_pair}. Since the diffusion fields in \Cref{eq:SO-Brownian-SDE} are linear and \(\sum_{a<b}J_{ab}^2=-(d-1)I_d\), the Stratonovich-to-It\^o conversion gives
    \[
        \begin{aligned}
            \mathrm d\mathcal    R_s & =\frac12\sum_{a<b}J_{ab}^2\mathcal R_s\,\mathrm ds+\sum_{a<b}J_{ab}\mathcal R_s\,\mathrm dB_s^{ab} \\
                                     & =-\frac{d-1}{2}\mathcal R_s\,\mathrm ds+\sum_{a<b}J_{ab}\mathcal R_s\,\mathrm dB_s^{ab}.
        \end{aligned}
    \]
    Consequently, one has for another standard Brownian motion $\tilde{B}_s^{ab}$ that
    \[
        \mathrm d\mathcal R_{8\Lambda s}=-4\Lambda(d-1)\mathcal R_{8\Lambda s}\,\mathrm ds+\sqrt{8\Lambda}\sum_{a<b}J_{ab}\mathcal R_{8\Lambda s}\,\mathrm d\tilde{B}_s^{ab}.
    \]
    Set \(\Omega_s:=X_s-X_s^*=\mathcal R_{8\Lambda s}\Omega_0\).  Noting that the center \(C\) is constant in \(s\) and \(X_s=C+\Omega_s/2\) and \(X_s^*=C-\Omega_s/2\). Then, one has
    \[
        \begin{aligned}
            \mathrm dX_s   & =\frac12\,\mathrm d\Omega_s=-2\Lambda(d-1)\Omega_s\,\mathrm ds+\sqrt{2\Lambda}\sum_{a<b}J_{ab}\Omega_s\,\mathrm d\tilde{B}_s^{ab},\qquad
            \mathrm dX_s^* & =-\frac12\,\mathrm d\Omega_s=-\mathrm dX_s.
        \end{aligned}
    \]
    Hence, the pair \(Y_s\) has continuous paths and satisfies
    \[
        |X_s|^2+|X_s^*|^2=|X_0|^2+|X_0^*|^2,
    \]
    so that the law $\nu_s\in \cP_2(\R^{2d})$.  Using the continuity, one has
    \[
        W_2^2(\nu_s,\nu_t)\leq\mathbb E\bigl[|X_s-X_t|^2+|X_s^*-X_t^*|^2\bigr]\longrightarrow0.
    \]
    Thus \(\nu_\cdot\in C([0,\tau];\mathcal P_2(\mathbb R^{2d}))\).

    Collecting these two equations yields the joint It\^o SDE
    \[
        \mathrm \dd Y_s=b(Y_s)\,\mathrm ds+\sum_{a<b}\sigma_{ab}(Y_s)\,\mathrm d\tilde{B}_s^{ab},
    \]
    where, for \(z=(v,v_*)\) and \(\omega=v-v_*\),
    \[
        b(z):=2\Lambda(d-1)\begin{pmatrix}-\omega\\ \omega\end{pmatrix},\qquad\sigma_{ab}(z):=\sqrt{2\Lambda}\begin{pmatrix}J_{ab}\omega\\-J_{ab}\omega\end{pmatrix}.
    \]
    Fix \(\varphi\in C_c^\infty(\mathbb R^{2d})\).  The multidimensional It\^o formula \cite[Theorem~18.18]{Kallenberg2021}  gives
    \[
        \mathrm d\varphi(Y_s)=\left[b\cdot\nabla\varphi+\frac12\sum_{a<b}\sigma_{ab}^{\mathsf T}D^2\varphi\,\sigma_{ab}\right](Y_s)\,\mathrm ds+\sum_{a<b}\nabla\varphi(Y_s)\cdot\sigma_{ab}(Y_s)\,\mathrm d\tilde{B}_s^{ab}.
    \]
    Here \(\nabla\varphi\) and \(D^2\varphi\) are taken in the full variable \(z=(v,v_*)\). We now identify the expression in brackets. In the coordinates \Cref{eq:center_relvel}, differentiation in \(\omega\) keeps \(c\) fixed, and
    \[
        \nabla_\omega=\frac12(\nabla_v-\nabla_{v_*}).
    \]
    For the Maxwell matrix \(A(\omega)=\Lambda(|\omega|^2I_d-\omega\otimes\omega)\), we have
    \[
        \sum_{a<b}(J_{ab}\omega)(J_{ab}\omega)^{\mathsf T}
        =|\omega|^2I_d-\omega\otimes\omega
        =\Lambda^{-1}A(\omega),
        \qquad
        \nabla_\omega\cdot A(\omega)=-\Lambda(d-1)\omega.
    \]
    Defining $\tilde{\varphi}(e; c, r):=\varphi(c+\frac{r}{2}e,
        c-\frac{r}{2}e )$. Substituting the drift and diffusion coefficients into the expression in brackets therefore gives
    \begin{equation}\label{eq:pair-generator}
        \begin{aligned}
            b\cdot\nabla_z\varphi
            +\frac12\sum_{a<b}
            \sigma_{ab}^{\mathsf T}D_z^2\varphi\,\sigma_{ab}
             & =-4\Lambda(d-1)\omega\cdot\nabla_\omega\varphi
            +4\sum_{p,q=1}^d A_{pq}(\omega)
            \frac{\partial^2\varphi}{\partial\omega_p\partial\omega_q} \\
             & =4\nabla_\omega\cdot
            \bigl(A(\omega)\nabla_\omega\varphi\bigr)= 4\Lambda \Delta_{\mathbb{S}^{d-1}} \tilde{\varphi}.
        \end{aligned}
    \end{equation}
    This clearly agrees with the significance of $\mathcal R$.
    Since \(A\) is symmetric, this operator is self-adjoint and thus the adjoint of the forward operator \(\mathcal L^*\) in \Cref{eq:pair_initial_value_problem}. Thus the drift term in It\^o's formula is \(\mathcal L\varphi(Y_s)\). The calculation also holds at \(\omega=0\), where both sides vanish. Moreover, \(\mathcal R_0=I_d\) gives \(\nu_0=f_k^{\otimes2}\).

    Since \(\mathcal L\varphi\) and \(\nabla\varphi\cdot\sigma_{ab}\) are bounded, the stochastic integral has zero expectation, and Fubini's theorem gives
    \[
        \mathbb E[\varphi(X_s,X_s^*)]-\mathbb E[\varphi(X_0,X_0^*)]=\int_0^s\mathbb E[\mathcal L\varphi(X_u,X_u^*)]\,\mathrm du,
    \]
    or equivalently,
    \begin{equation}\label{eq:joint_law_weak_equation}
        \int_{\mathbb R^{2d}}\varphi\,\mathrm d\nu_s-\int_{\mathbb R^{2d}}\varphi\,\mathrm d(f_k^{\otimes2})=\int_0^s\int_{\mathbb R^{2d}}\mathcal L\varphi\,\mathrm d\nu_u\,\mathrm du.
    \end{equation}
    Since this identity holds for every \(\varphi\in C_c^\infty(\mathbb R^{2d})\), then $\nu_{\cdot}$ is a solution to the equation in the weak sense
    \[
        \partial_s\nu_s=\mathcal L^*\nu_s,\qquad \nu_0=f_k^{\otimes2}.
    \]

    By the uniqueness in \(C([0,\tau];\mathcal P_2(\mathbb R^{2d}))\), equipped with the \(W_2\) topology in \Cref{lem:pair-forward-uniqueness}, \(\nu_s=F_{t_k+s}\) for \(0\le s<\tau\), and continuity at the terminal time gives \(\nu_\tau=F_{k+1}^-\).

    Since \(P_\tau=\mathcal R_{8\Lambda \tau}\) and
    \[
        X_\tau=\frac{X_0+X_0^*}{2}+\frac12P_\tau(X_0-X_0^*)=A_{P_\tau}X_0+B_{P_\tau}X_0^*,
    \]
    one thus has
    \[
        f_{k+1}=\pi_1F_{k+1}^-=\pi_1\nu_\tau=\Law(A_{P_\tau}X_0+B_{P_\tau}X_0^*)=S_\tau f_k,
    \]
    which proves \Cref{eq:rho_recursion}.

    The pair map conserves pair momentum and energy.  Symmetry of the two outputs therefore gives
    \begin{equation}
        M_2(f_k):=\int_{\R^d}|v|^2f_k(\dd v)=M_2(f_0)\qquad\text{for every }k.\label{eq:M2-conservation}
    \end{equation}
    We have then verified the claims.
\end{proof}

\begin{lemma}[Finite-time propagation of a higher moment]\label[lemma]{lem:moment}
    Let $q>2$ and $M_q(f_0)<\infty$.  Then there is a constant $C=C(d,q,\Lambda)$ such that
    \begin{equation}
        \sup_{k\tau\leq T}M_q(f_k)\leq e^{CT}M_q(f_0).\label{eq:moment-bound}
    \end{equation}
    In particular, the constant is independent of $\tau$.
\end{lemma}

\begin{proof}
    Fix $k$ and suppose that $M_q(f_k)<\infty$. We first estimate the collision increment for fixed pre-collisional velocities $v,v_*$. Write $c=(v+v_*)/2$, $r=|v-v_*|$, and, when $r>0$, $e_0=(v-v_*)/r$; if $r=0$, choose any $e_0\in\mathbb S^{d-1}$. Define
    \[
        f_{c,r}(e)
        :=|c+(r/2)e|^q+|c-(r/2)e|^q.
    \]
    Clearly, $f_{c,r}(e_0)=|v|^q+|v_*|^q$. If $r=0$, then $f_{c,0}$ is constant and independent of $e_0$.

    Let $e_s:=\mathcal R_{8\Lambda s}e_0$, $0\leq s\leq \tau$, so that $c\pm(r/2)e_s$ are the two velocities in \Cref{eq:candidate_pair} for this initial pair. The process $e_s$ has generator $4\Lambda\Delta_{\mathbb S^{d-1}}$. Denote by $\mathbb E_{\mathrm{rot}}$ expectation over the rotation with $v,v_*$ fixed. Since $q>2$, the function $g(x)=|x|^q$ belongs to $C^2(\mathbb R^d)$. Then, It\^o's formula gives
    \[
        \mathbb E_{\mathrm{rot}}[f_{c,r}(e_\tau)] =f_{c,r}(e_0)+4\Lambda\int_0^\tau\mathbb E_{\mathrm{rot}}\left[\Delta_{\mathbb S^{d-1}}f_{c,r}(e_s)\right]\,\mathrm ds.
    \]
    For $H\in C^2(\mathbb R^d)$ and $e\in\mathbb S^{d-1}$,
    \[
        \begin{aligned}
            \Delta_{\mathbb S^{d-1}}(H|_{\mathbb S^{d-1}})(e)
             & =\Delta H(e)-e^{\mathsf T}D^2H(e)e-(d-1)e\cdot\nabla H(e)                        \\
             & =\operatorname{tr}\!\left[(I_d-e\otimes e)D^2H(e)\right]-(d-1)e\cdot\nabla H(e).
        \end{aligned}
    \]
    Setting \(\varrho=r/2\) and \(x_\pm=c\pm\varrho e\), we obtain
    \[
        \Delta_{\mathbb S^{d-1}}g(c\pm\varrho e)
        =\varrho^2\operatorname{tr}\!\left[(I_d-e\otimes e)D^2g(x_\pm)\right]
        \mp(d-1)\varrho\,e\cdot\nabla g(x_\pm).
    \]
    Since
    \[
        \begin{aligned}
            \left|\operatorname{tr}\!\left[(I_d-e\otimes e)D^2g(x)\right]\right|
             & =q(d-1)|x|^{q-2}+q(q-2)|x|^{q-4}\bigl(|x|^2-(e\cdot x)^2\bigr) \\
             & \leq q(d+q-3)|x|^{q-2},                                        \\
            |e\cdot\nabla g(x)|
             & =q|x|^{q-2}|e\cdot x|\leq q|x|^{q-1}.
        \end{aligned}
    \]
    The calculation is valid for \(x\neq0\), and the bounds extend to \(x=0\) by continuity. One then can bound easily that
    \begin{equation}
        |\Delta_{\mathbb S^{d-1}}f_{c,r}(e)|
        \leq\sum_{\pm}\left|\Delta_{\mathbb S^{d-1}}g(c\pm\varrho e)\right| \leq C_{d,q}(|v|^q+|v_*|^q).
        \label{eq:spherical-laplacian-moment}
    \end{equation}

    Now sample $v,v_*$ independently with law $f_k$, independently of the rotation, and write $\mathbb E_{v,v_*}$ for expectation over this initial pair. Set $v'=c+(r/2)e_\tau$ and $v_*'=c-(r/2)e_\tau$. By \Cref{eq:Sh,eq:rho_recursion} and symmetry under exchange of the initial velocities, both outputs have law $f_{k+1}$. Hence
    \[
        \begin{aligned}
            2M_q(f_{k+1}) =\mathbb E_{v,v_*}\left[\mathbb E_{\mathrm{rot}}[|v'|^q+|v_*'|^q]\right]\leq(1+C\tau)\mathbb E_{v,v_*}[|v|^q+|v_*|^q]   =2(1+C\tau)M_q(f_k).
        \end{aligned}
    \]
    Thus $M_q(f_{k+1})\leq(1+C\tau)M_q(f_k)$, and \Cref{eq:moment-bound} follows with a constant independent of $\tau$.
\end{proof}

\subsection{Wasserstein distances and relative entropy}

We will make use of the Wasserstein distances to gauge the closeness between the empirical measure of the particles and the
target probability distribution. Let $\cP(\Omega)$ be the set of all probability measures on a Euclidean space $\Omega$ and $\cP_p(\Omega)$ indexed by some $p\ge 1$ be the set of probability measures having finite $p$-moment, i.e.,
\begin{equation}\label{eq:P2}
    \mathcal{P}_{p}\left(\Omega\right):=\left\{\mu \in \mathcal{P}\left(\Omega\right): \int_{\Omega}\mathrm{dist}(x, 0)^p \mu(d x)<\infty\right\}.
\end{equation}
Clearly, $\cP_p\subset \cP_2$ for $p\ge 2$.
For $p\ge 1$ and any $\mu, \nu\in\mathcal{P}_{p}(\Omega)$, the $p$-Wasserstein distance is defined by
\begin{equation*}
    {W}_{p}(\mu, \nu):=\left(\inf _{\pi \in \Pi(\mu, \nu)} \int_{\Omega \times \Omega}\mathrm{dist}(x, y)^{p} \pi(d x, d y)\right)^{1/p},
\end{equation*}
where $\mathrm{dist}(x, y)$ is the distance between $x$ and $y$, and $\Pi(\mu, \nu)$ is the set of couplings of $\mu$ and $\nu$ (i.e., joint distributions with marginals to be $\mu$ and $\nu$ respectively). The space $\mathcal{P}_{p}\left(\Omega\right)$ is a Polish space under the $p$-Wasserstein metric\cite[Theorem~6.18]{Villani2009}.

In this work, we will choose mostly $\Omega=\R^d$ for the velocity space and
\[
    \mathrm{dist}(v_1, v_2)=|v_1-v_2|
\]
is chosen to be the Euclidean distance. In some proofs later, we may consider $\Omega=(\R^d)^{N}$
to be the space for the $N$-particle configuration.

We will also consider the relative entropy (KL divergence) in later sections to measure the discrepancy between two probability measures, which is defined by the following
\begin{gather}\label{eq:relativeentropy}
    \cH(\mu\mid \nu) = \left\{
    \begin{aligned}
         & \int_{\Omega} \log \frac{\mathrm{d} \mu}{\mathrm{d}\nu} \mathrm{d}\mu, & \text{if} ~ \mu \ll \nu, \\
         & \infty,                                                                & \text{else}.
    \end{aligned}\right.
\end{gather}
Here, $\frac{\mathrm{d} \mu}{\mathrm{d}\nu}$ denotes the Radon-Nikodym derivative of $\mu$ with respect to $\nu$.
The relative entropy is a non-negative quantity by Jensen's inequality and achieves zero only if $\mu= \nu$. Moreover, the relative entropy is jointly convex in its arguments and is lower semicontinuous with respect to weak convergence of probability measures\cite[Theorems~11 and~19]{vanErvenHarremoes2014}.

A class of important functional inequalities related to relative entropy is the transportation inequalities \cite{talagrand1996transportation,bobkov2001hypercontractivity}. If $\nu$ satisfies a log-Sobolev inequality, then one can control the $W_2$ distance by the relative entropy as follows
\cite[Corollary~3.1]{bobkov2001hypercontractivity}:
\begin{equation}
    W_1(\mu, \nu)\le W_2(\mu,\nu) \leq C\sqrt{ \cH(\mu |\nu)},
\end{equation}
where $C$ depends only on the log-Sobolev constant.
Besides, if $\nu$ is sub-Gaussian, then one has by the weighted Csisz\'ar-Kullback-Pinsker inequality \cite[Corollary~2.4]{bolley2005weighted}, one has
\begin{equation}
    W_1(\mu,\nu) \leq C_{\nu}\sqrt{ \cH(\mu |\nu)},
\end{equation}
where
\begin{equation}\label{eq:W1condition}
    C_\nu:=2 \inf _{\alpha>0}\left(\frac{1}{2 \alpha}\left(1+\log \int_{\Omega} e^{\alpha|x|^2} d \nu(x)\right)\right)^{\frac{1}{2}}<+\infty.
\end{equation}

\subsection{Main results}

Our main results consist of two parts. The first is the large particle limit which shows that $f^N_k$ defined in \Cref{eq:fkNdef}
converges in Wasserstein distance to $f_k$. The second is concerned with the convergence of $f_k$ to the solution of Landau equation. To state the results, we introduce for \(d\ge2\) and \(N\ge2\) the dimension-dependent rate
\begin{equation}\label{eq:dimension-rate}
    \varepsilon_d(N):=
    \begin{cases}
        N^{-1/2},           & 2\le d<4, \\
        N^{-1/2}[1+\log N], & d=4,      \\
        N^{-2/d},           & d>4.
    \end{cases}
\end{equation}
The first main result is then given as follows.
\begin{theorem} \label{thm:main}
    Let $d\ge2$ be an integer, $\gamma=0$, $T>0$, $\Lambda>0$, $q>4$, and $f_0\in\mathcal P_q(\R^d)$. For every integer $K\geq1$, set $\tau=T/K$. For every even $N\geq2$, initialize $(V_0^i)_{i=1}^N$ independently with law $f_0$. There exists a finite constant $C=C(T,d,q,\Lambda)$ such that
    \begin{equation}
        \max_{0\leq k\leq K}
        \E\Wtwo^2(f_k^N, f_k)
        \leq C\,M_q(f_0)^{2/q}\varepsilon_d(N).
        \label{eq:main-uniform-bound}
    \end{equation}
    In particular, $C$ is independent of $N$, $K$, and $\tau$.
\end{theorem}
 
 The second main result is the entropic convergence of $f_k$ to the solution $f$ of the Landau equation \Cref{eq:landau}.
Denote
 \[
 \langle v\rangle=(1+|v|^2)^{1/2}.
 \]
  For the entropic convergence estimate, we impose the following stronger assumptions on the common initial density $f_0$.
\begin{assumption}\label{ass:step2-initial}
    The scheme and the Landau equation start from the same strictly positive $C^3$ probability density $f_0$ on $\R^d$, for a fixed $d\ge2$.  For some $m_0\ge1$
    and $K_0<\infty$,
    \[
        |\nabla_v\log f_0(v)|+|\nabla_v^2\log f_0(v)|
        \le K_0\langle v\rangle^{m_0}.
    \]
    For some $q>6$, with $m=m_0(m_0+1)$ and $p_*=\max\{16,8q/(q-6),2m+4\}$, one has
    \[
        \sum_{j=1}^3\int f_0\left|\frac{\nabla^j f_0}{f_0}\right|^{q/j}\dd v <\infty,\quad M_{p_*}(f_0)<\infty,
    \]
    where for a function $h$
    \[
        |\nabla^j h|^2=\sum_{i_1,\ldots,i_j=1}^d |\partial_{v_{i_1}}\cdots \partial_{v_{i_j}}h|^2=\sum_{|\alpha|=j}\frac{j!}{\alpha!}
        |\partial^{\alpha}h|^2.
    \]
\end{assumption}

The second main result is given by the following.
\begin{theorem}[Entropic convergence for Maxwell molecules]
    \label{thm:step2-entropy}
    Under \Cref{ass:step2-initial}, let $f(t,v)$ solve \Cref{eq:landau}
    with $\gamma=0$ and initial density $f_0$. Then, one has 
    \begin{equation}\label{eq:step2-main}
        \sup_{k: k\tau \le T} \cH(f_k \mid f(t_k,\cdot)) \le C_T\tau^2,
    \end{equation}
    where the constant $C_T$ is independent of $\tau$ and depends only on $T$, $d$,
    $\Lambda$, and the initial quantities in \Cref{ass:step2-initial}.
\end{theorem}

Since $p_*\ge16>4$, \Cref{thm:main} applies with moment order $p_*$ under \Cref{ass:step2-initial}. Together with \Cref{thm:step2-entropy}, one has the following overall estimate for the particle scheme.
\begin{corollary}\label{cor:errorinW}
    Assume \Cref{ass:step2-initial}. Let $N\ge2$ be even, $K\ge1$, and $\tau=T/K$, and initialize the particles independently with common law $f_0$. The following estimates hold for every $0\le k\le K$.

    The particle system satisfies
    \[
        \E W_1^2(f_k^N, f(t_k,\cdot))\le C_T\tau^2+C_T \varepsilon_d(N).
    \]
    If $f(t,\cdot)$ satisfies a uniform log-Sobolev inequality, then one has
    \[
        \E W_2^2(f_k^N, f(t_k,\cdot))\le C_T\tau^2+C_T \varepsilon_d(N).
    \]
    Here $C_T$ is independent of $N$, $K$, and $\tau$; it depends on the initial quantities in \Cref{ass:step2-initial}, $T$, $d$, $\Lambda$, and, for the $W_2$ estimate, the uniform log-Sobolev constant.
      If $d=3$, the bound (for $p=1$ or $p=2$) reads
    \[
        \E  W_p^2(f_k^N, f(t_k,\cdot))\le C_T(\tau^2+N^{-1/2}).
    \]

\end{corollary}

\begin{figure}[H]
    \centering
    \includegraphics[width=0.9\textwidth]{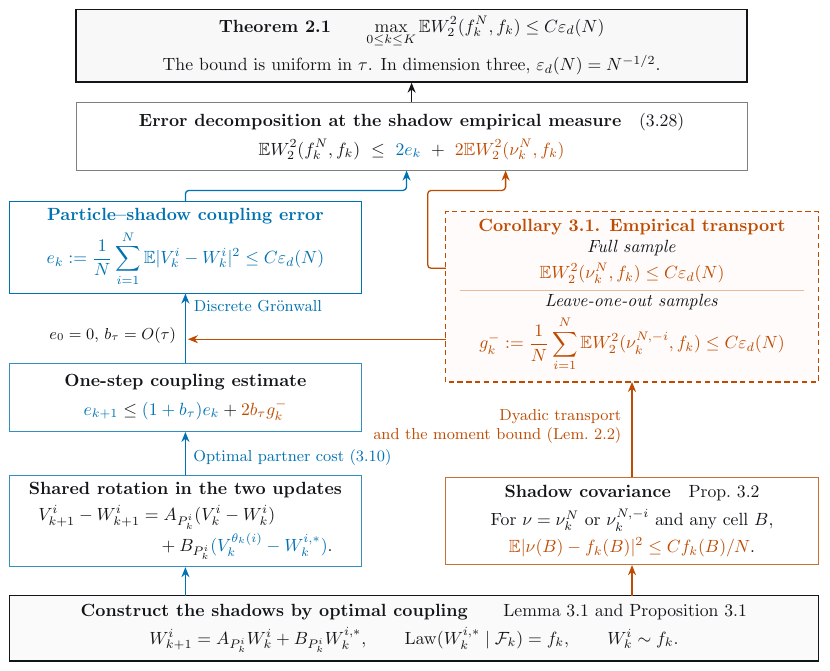}
    \caption{Proof of \Cref{thm:main}. The full-sample estimate enters the error decomposition directly, and the leave-one-out estimate enters the coupling recursion. Here $\mathcal F_k$ is the history before the matching at step $k$, and $C$ is independent of $N$ and $\tau$.}
    \label{fig:proof-structure}
\end{figure}

\begin{figure}[H]
    \centering
    \includegraphics[width=0.9\textwidth]{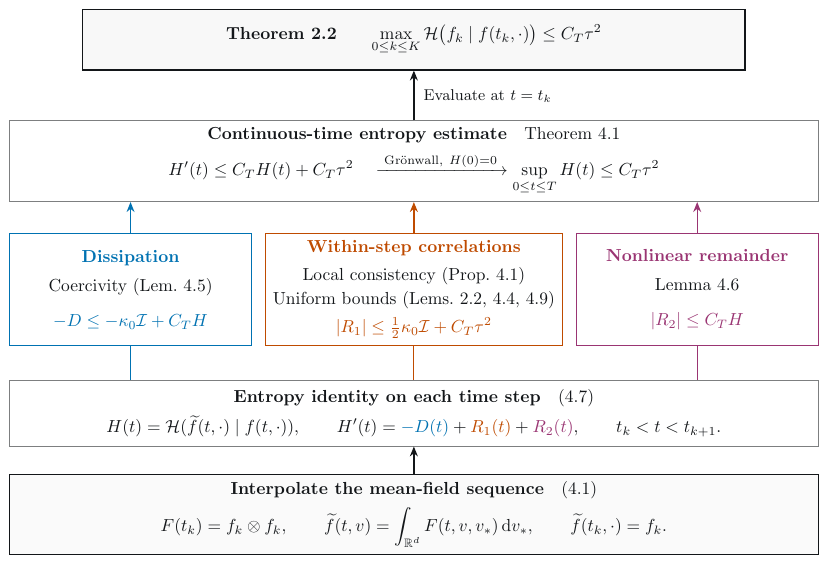}
    \caption{Proof of \Cref{thm:step2-entropy}. Here $\cI$ denotes the relative Fisher information of $\widetilde f$ with respect to $f$, and $\kappa_0>0$ is the coercivity constant. The dissipation absorbs the $\cI$ term in the correlation estimate, leaving a forcing term of order $\tau^2$.}
    \label{fig:step2-proof}
\end{figure}

Below, we give the outlines of the proofs for the main results.
\Cref{sec:large_particle} compares the particle system with the mean-field law $f_k$ by coupling, as outlined in \Cref{fig:proof-structure}. \Cref{sec:entropic} compares the mean-field interpolation with the Landau solution by relative entropy, as outlined in \Cref{fig:step2-proof}.

Adapting the construction of \cite{BorghiPareschi2025}, we couple each particle $V_k^i$ to a shadow $W_k^i$ with initial velocity \(W_0^i=V_0^i\). Using the Maxwell representation of \Cref{subsec:sbmbasics,subsec:limitdynamics}, we evolve both with the same rotation \(P_k^i\),
\begin{equation}\label{eq:original_shadow}
    \begin{aligned}
        V_{k+1}^i
         & =A_{P_k^i}V_k^i+B_{P_k^i}V_k^{\theta_k(i)}, \\
        W_{k+1}^i
         & =A_{P_k^i}W_k^i+B_{P_k^i}W_k^{i,*}.
    \end{aligned}
\end{equation}
By \Cref{lem:labelled-OT}, the conditional law of \((W_k^{i,*},V_k^{\theta_k(i)})\) given the pre-step history is an optimal quadratic-cost coupling of \(f_k\) and \(f_k^{N,-i}\), the empirical law of the other \(N-1\) particle velocities.
Each \(W_k^i\) has law \(f_k\) by \Cref{prop:correct-marginal}, although the matching and shared rotations can induce dependence between shadows.

Let \(\nu_k^N\) be the empirical law of the shadows and \(e_k=N^{-1}\sum_{i=1}^N\mathbb E|V_k^i-W_k^i|^2\) the mean-square coupling error. The labelwise coupling and the triangle inequality give
\[
    \mathbb E W_2^2(f_k^N, f_k)
    \leq 2e_k+2\mathbb E W_2^2(\nu_k^N, f_k).
\]
The shared rotations and optimal partner couplings in \Cref{eq:original_shadow} yield
\[
    e_{k+1}\leq(1+b_\tau)e_k+2b_\tau g_k^-,
    \qquad e_0=0,
    \qquad 0\leq b_\tau\leq C\tau,
\]
where \(g_k^-\) is the average of \(\mathbb E W_2^2(\nu_k^{N,-i},f_k)\) over \(i\), and \(\nu_k^{N,-i}\) is the empirical law with shadow \(i\) omitted.

For bounded Borel functions \(\phi,\psi:\R^d\to\R\), \Cref{prop:shadow-covariance} gives
\[
    \left|\operatorname{Cov}\bigl(\phi(W_k^i),\psi(W_k^j)\bigr)\right|
    \leq\frac5{N-1}\sqrt{\operatorname{Var}_{f_k}(\phi)\operatorname{Var}_{f_k}(\psi)},
    \qquad i\ne j.
\]
For a cell indicator \(\phi=\psi=\one_B\), the variance factor is \(f_k(B)(1-f_k(B))\). The dyadic transport estimate of \cite[Lemmas~5--6]{FournierGuillin2015} and the moment bound of \Cref{lem:moment} then give the two bounds of order \(\varepsilon_d(N)\) in \Cref{cor:shadow-empirical}. The full-sample bound controls \(\mathbb E W_2^2(\nu_k^N,f_k)\) directly. The leave-one-out bound controls \(g_k^-\), so \(b_\tau=O(\tau)\) and discrete Gr\"onwall give the same rate for \(e_k\), uniformly for \(k\tau\le T\). The error decomposition now yields \Cref{thm:main}.

For the entropic convergence of $f_k$ to the Landau solution $f$, we  differentiate the relative entropy $H(t)=\cH(\widetilde f(t,\cdot)\mid f(t,\cdot))$ within each step $t\in (t_k, t_{k+1})$, where $\tilde{f}$ is the natural interpolation as the marginal of the two particle evolution in \eqref{eq:pair_initial_value_problem}:
\[
    H'(t)=-D(t)+R_1(t)+R_2(t), \qquad t_k<t<t_{k+1}.
\]
By \Cref{lem:relativecoercivity}, the dissipation $D$ controls the relative Fisher information up to a multiple of $H$. The local consistency estimate of \Cref{prop:Theta}, together with the moment and derivative bounds of \Cref{lem:moment,lem:step2-propagation,lem:Bbound}, controls $R_1$, the error arising from the within-step difference $F-\widetilde f\otimes\widetilde f$. \Cref{lem:prooflem2} bounds the nonlinear remainder $R_2$ by a multiple of $H$. Since $H(0)=0$, Gr\"onwall's inequality gives $H(t)\le C_T\tau^2$. At the grid times this gives \Cref{thm:step2-entropy}. Converting the entropy estimate to a Wasserstein bound and applying \Cref{thm:main} yields \Cref{cor:errorinW}.

\section{Proof of the large particle limit}\label{sec:large_particle}

In this section, we begin to prove \Cref{thm:main} in detail.

\subsection{The shadow system}
We adapt the coupling strategy of \cite{BorghiPareschi2025} by constructing the shadow system outlined above, \(\boldsymbol W_k=(W_k^1,\ldots,W_k^N)\). At each step, \(W_k^i\) collides with an auxiliary partner \(W_k^{i,*}\), using the same random rotation as the corresponding finite particle \(V_k^i\). Conditional on the current configuration \(\boldsymbol V_k=(V_k^1,\ldots,V_k^N)\), the pair \((W_k^{i,*},V_k^{\theta_k(i)})\) realizes an optimal quadratic-cost coupling between \( f_k\) and the leave-one-out empirical measure \(f_k^{N,-i}\). The resulting shadow particles have the correct one-particle marginals: \(\operatorname{Law}(W_k^i)=f_k\) for every \(i\) and \(k\).

An additional difficulty arises from the perfect-matching structure. In the Nanbu setting of \cite{BorghiPareschi2025}, independent partner selections and collision randomizations allow the nonlinear comparison particles to remain i.i.d. In our scheme, the reciprocal matching and the rotation shared by each matched pair introduce dependence between the shadow particles. Thus, although each \(W_k^i\) has law \(f_k\), their joint law need not be \(f_k^{\otimes N}\). We estimate this dependence through the covariance of shadow observables and use the resulting bound to control the shadow empirical-measure error.

We begin by constructing the probability kernels used to sample the shadow partners \(W_k^{i,*}\).
For the perfect matching at any step $k$, conditionally on the past, $\theta_k(i)$ is uniform on the $N-1$ labels other than $i$. Hence we consider the leave-one-out empirical measure: for a deterministic configuation $\boldsymbol{v}=(v^1,\cdots,v^N)$ and a label $i$, set
\[
    \mu_{\boldsymbol{v}}^{-i}:=\frac{1}{N-1}\sum_{\ell\neq i}\delta_{v^\ell}.
\]

The following lemma gives the probability kernels \(G_k^{i,\ell}(\boldsymbol v,\cdot)\in\mathcal P(\mathbb R^d)\), representing the conditional law of an auxiliary partner velocity for particle \(i\), given the configuration \(\boldsymbol v\) and the selected partner label \(\ell\neq i\). More precisely, if \(\ell\) is sampled uniformly from the labels other than \(i\), and \(w\) is then sampled from \(G_k^{i,\ell}(\boldsymbol v,\cdot)\), the resulting pair \((v^\ell,w)\) realizes an optimal coupling between \(\mu_{\boldsymbol v}^{-i}\) and \( f_k\).
\begin{lemma}[Labelled optimal-coupling kernels]
    \label[lemma]{lem:labelled-OT}
    For every $k,i$ and configuration $\bm v$, one can choose probability kernels $G_k^{i,\ell}(\bm v,\dd w)$, $\ell\neq i$, such that
    \begin{align}
        \frac1{N-1}\sum_{\ell\neq i} G_k^{i,\ell}(\bm v,\dd w)                         & =f_k(\dd w), \label{eq:G-marginal}                 \\
        \frac1{N-1}\sum_{\ell\neq i} \int_{\R^d}|v^\ell-w|^2 G_k^{i,\ell}(\bm v,\dd w) & =\Wtwo^2(\mu_{\bm v}^{-i}, f_k). \label{eq:G-cost}
    \end{align}
    The kernels may be chosen jointly measurable in $\bm v$ and strictly permutation-equivariant: for every permutation $\sigma$,
    \begin{equation}
        G_k^{\sigma(i),\sigma(\ell)}
        (\sigma\bm v,\cdot)
        =G_k^{i,\ell}(\bm v,\cdot).
        \label{eq:G-equivariance}
    \end{equation}
\end{lemma}
We follow the construction in
\cite[Lemma~4.1]{BorghiPareschi2025} with $p=2$, adapted to the leave-one-out empirical measure. We briefly recall the argument and verify the permutation-equivariance property.
In fact, the construction is based on conditioning an optimal coupling. Suppose $\pi=\Phi_k(\mu_{\boldsymbol{v}}^{-i})$ is an optimal transport plan from $f_k$ to $\mu_{\boldsymbol{v}}^{-i}$. Let \((W,Z)\sim\pi\). Then \(W\sim f_k\) and \(Z\sim\mu_{\boldsymbol v}^{-i}\). The kernel \(G_k^{i,\ell}(\bm v,\cdot)\) represents the conditional law of \(W\) given \(Z=v^\ell\). This suggests the conditional-probability formula:
\begin{align} \label{eq:conditional_G}
    G_k^{i,\ell}(\boldsymbol{v}, A) = \frac{\Phi_k(\mu_{\boldsymbol{v}}^{-i})(A\times \{v^\ell\})}{\mu_{\boldsymbol{v}}^{-i}(\{v^\ell\})}.
\end{align}
By \cite[Corollary~5.22]{Villani2009}, there exists a Borel selector of optimal transport plans
\begin{align}\label{eq:Borel-OT-selector}
    \Phi_k:\mathcal P_2(\R^d)\longrightarrow\mathcal P_2(\R^{2d}), \qquad \Phi_k(\mu)\in \Pi_{\operatorname{opt}}(f_k,\mu).
\end{align}
In \Cref{eq:conditional_G}, the numerator is Borel in \(\bm v\) by parameter-dependent kernel integration, while the denominator is Borel and at least \(1/(N-1)\); hence \(G_k^{i,\ell}\) is a Borel probability kernel. The marginal and optimal-cost identities then follow by reconstructing the selected optimal plan from these kernels, while permutation equivariance follows because the construction depends only on \(\mu_{\boldsymbol v}^{-i}\) and \(v^\ell\).

By the kernel representation lemma \cite[Lemma~4.22]{Kallenberg2021}, we have a Borel map \(\mathfrak{R}:\mathcal{P}(\mathbb{R}^d)\times [0,1]\to\R^d\) such that for all \(\xi\sim \operatorname{U}(0,1)\), we have
\[
    \mathrm{Law}(\mathfrak{R}(\eta,\xi)) =\eta.
\]
For fixed \(k, i, \ell\), define
\[
    g_k^{i,\ell}(\boldsymbol{v},u):=\mathfrak{R}(G_k^{i,\ell}(\boldsymbol{v},\cdot), u).
\]
By \Cref{lem:labelled-OT}, $G_k^{i,\ell}$ is a Borel probability kernel, hence $g_k^{i,\ell}$ is jointly Borel in $(\boldsymbol{v}, u)$. So for every bounded Borel function $\phi$, we have
\begin{align}
    \int_{[0,1]}\phi(g^{i,\ell}_k(\boldsymbol{v}, u))\dd u = \int_{\R^d}\phi(\omega) G_k^{i,\ell}(\boldsymbol{v},\dd \omega).
\end{align}

With this sampler, we can now define shadow collision partners and further define the shadow system $\boldsymbol{W}_k$.
Let
\[
    \mathcal F_k
    :=\sigma\bigl(\boldsymbol V_0,\theta_r,P_r^i,U_r^i:
    0\le r<k,\ 1\le i\le N\bigr)
\]
be the pre-step \(\sigma\)-field. At step \(k\), sample a uniform perfect matching \(\theta_k\) independently of \(\mathcal F_k\). Conditionally on \(\mathcal F_k\vee\sigma(\theta_k)\), the rotations on distinct matching edges are independent with common law \(\Law(P_\tau)\).
Also, denote
\[
    \tilde{\cF}_k=\cF_k\vee\sigma(\theta_k,P_k^1,\ldots,P_k^N).
\]
Then, $\boldsymbol{W}_k$ is constructed as follows.
\begin{itemize}
    \item Set $\boldsymbol{W}_0=\boldsymbol{V}_0\sim f_0^{\otimes N}$.
    \item Assume that we have defined $\boldsymbol{W}_0,\cdots ,\boldsymbol{W}_k$. Choose mutually independent $U_k^1,\cdots, U_k^N\sim \operatorname{U}(0,1)$ such that the entire vector is independent of the joint $\sigma$-field $\tilde{\cF}_k$.
          Define
          \begin{align}
              W_k^{i,*} := g_k^{i,\theta_k(i)}(\boldsymbol{V}_k,U_k^i).
          \end{align}

    \item Then, $\boldsymbol{W}_{k+1}$ is obtained by the second equation in \Cref{eq:original_shadow}.
\end{itemize}

Note that $W_k^{i,*} $ is a random variable because $g_k^{i,\ell}$ is jointly Borel and the label $\theta_k(i)$ takes only finitely many values.  Clearly, $W_k^{i,*}$ need not coincide with any existing shadow velocity $W_k^\ell$.
By the construction and applying \Cref{eq:G-marginal}, one has
\begin{equation}
    \begin{split}
        \E\left[\phi(W_k^{i,*})\mid\cF_k\right]
         & =\frac1{N-1}\sum_{\ell\neq i}
        \int_{\R^d}\phi(w)G_k^{i,\ell}(\bm V_k,\dd w) \\
         & =\int_{\R^d}\phi(w) f_k(\dd w).
        \label{eq:partner-conditional-expectation}
    \end{split}
\end{equation}
Consequently,
\begin{equation}
    \Law(W_k^{i,*}\mid\cF_k)=f_k.
    \label{eq:partner-law}
\end{equation}
The measure \(f_k\) is deterministic and thus \(W_k^{i,*}\) is independent of \(\mathcal F_k\), and in particular independent of  \(W_k^i\).
Similarly,  the construction and \Cref{eq:G-cost} gives
\begin{equation}
    \E\left[|V_k^{\theta_k(i)}-W_k^{i,*}|^2\mid\cF_k\right]=\Wtwo^2(\mu_{\bm V_k}^{-i}, f_k)
    =\Wtwo^2(f_k^{N,-i}, f_k),
    \label{eq:partner-optimal-cost}
\end{equation}
where
\begin{equation}
    f_k^{N,-i}:=\frac1{N-1}\sum_{\ell\neq i}\delta_{V_k^\ell}.
    \label{eq:finite-loo}
\end{equation}

By the independence between $W_k^{i,*}$ and $W_k^i$ and the symmetry, one can easily obtain the following by induction whose proof is omitted. It is worth noting that, however, these shadow vectors $\{W_k^i\}_{i=1}^N$ are not independent in general.
\begin{proposition}[Correct one-particle marginal]
    \label[proposition]{prop:correct-marginal}
    For every $i$ and $k$,
    \begin{equation}
        W_k^i\sim  f_k.
        \label{eq:correct-marginal}
    \end{equation}
    For each fixed $k$, the shadow vector $\bm W_k=(W_k^1,\ldots,W_k^N)$ is exchangeable: for every permutation $\sigma$ of $\{1,\ldots,N\}$,
    \[
        \Law(\sigma\boldsymbol{W}_k)=\Law\bigl(W_k^{\sigma^{-1}(1)},\ldots,W_k^{\sigma^{-1}(N)}\bigr)=\Law(W_k^1,\ldots,W_k^N)=\Law(\boldsymbol{W}_k).
    \]
\end{proposition}

\subsection{Covariance and empirical transport for the shadow system}

We first control the covariance of bounded shadow observables. The resulting estimate will be applied to cell indicators to bound the full and leave-one-out shadow empirical measures.

\begin{proposition}[Shadow covariance]\label[proposition]{prop:shadow-covariance}
    For distinct labels $i,j$, every step $k$, and bounded Borel functions $\phi,\psi:\R^d\to\R$, one has
    \begin{equation}\label{eq:shadow-covariance}
        \left|\operatorname{Cov}\bigl(\phi(W_k^i),\psi(W_k^j)\bigr)\right|
        \leq\frac5{N-1}
        \sqrt{\operatorname{Var}_{f_k}(\phi)\operatorname{Var}_{f_k}(\psi)}.
    \end{equation}
\end{proposition}

\begin{proof}
    Fix a terminal step $k$, distinct labels $i,j$, and bounded Borel functions $\phi,\psi$. For a bounded Borel function $h$, define the one-step transition operator
    \[
        (Q_rh)(x):=\E h(A_Px+B_PZ),
        \qquad Z\sim f_r,\quad P\sim\Law(P_\tau),\quad P\perp Z.
    \]
    Let $P_{k,k}$ be the identity and set $P_{r,k}=Q_rP_{r+1,k}$ for $r<k$.
    For every label $\ell$, \cref{eq:partner-law} and the sampling construction imply that, conditionally on $\mathcal F_r$, $P_r^\ell$ and $W_r^{\ell,*}$ are independent with respective laws $\Law(P_\tau)$ and $f_r$. Thus \cref{eq:original_shadow} yields
    \begin{equation}\label{eq:shadow-transition-identity}
        \E[h(W_{r+1}^\ell)\mid\mathcal F_r]=(Q_rh)(W_r^\ell).
    \end{equation}
    Define processes
    \[
        \begin{aligned}
            M_r^i & :=(P_{r,k}\phi)(W_r^i), \qquad L_r^j & :=(P_{r,k}\psi)(W_r^j).
        \end{aligned}
    \]
    By \cref{eq:shadow-transition-identity}, for every $r<k$,
    \[
        \begin{aligned}
            \E[M_{r+1}^i\mid\mathcal F_r]
             & =\E[(P_{r+1,k}\phi)(W_{r+1}^i)\mid\mathcal F_r] \\
             & =(Q_rP_{r+1,k}\phi)(W_r^i)=M_r^i,
        \end{aligned}
    \]
    and similarly $\E[L_{r+1}^j\mid\mathcal F_r]=L_r^j$. Both processes are adapted and bounded, so they are martingales in the joint filtration $(\mathcal F_r)_{0\le r\le k}$. Since $P_{k,k}$ is the identity, $M_k^i=\phi(W_k^i)$ and $L_k^j=\psi(W_k^j)$.
    Their initial values are independent, since $W_0^i,W_0^j$ are independent and the transition operators are deterministic. Writing
    \[
        \Delta M_r^i:=M_{r+1}^i-M_r^i,\qquad
        \Delta L_r^j:=L_{r+1}^j-L_r^j,
    \]
    martingale orthogonality and the independence of the initial values give
    \begin{equation}\label{eq:shadow-covariance-increments}
        \begin{aligned}
            \operatorname{Cov}\bigl(\phi(W_k^i),\psi(W_k^j)\bigr) & =\operatorname{Cov}\bigl(M_k^i,L_k^j\bigr)       \\
                                                                  & =\operatorname{Cov}\bigl(M_0^i,L_0^j\bigr)
            +\sum_{r=0}^{k-1}\E[\Delta M_r^i\,\Delta L_r^j]                                                          \\
                                                                  & =\sum_{r=0}^{k-1}\E[\Delta M_r^i\,\Delta L_r^j].
        \end{aligned}
    \end{equation}

    We next show that, for each $r<k$,
    \begin{equation}\label{eq:shadow-increment-covariance}
        \left|\E[\Delta M_r^i\,\Delta L_r^j]\right|
        \leq\frac5{N-1}\|\Delta M_r^i\|_{L^2}\|\Delta L_r^j\|_{L^2},
    \end{equation}
    where $\|Z\|_{L^2}^2=\E|Z|^2$. Fix $r<k$. We split the expectation according to whether $i$ and $j$ belong to the same matching edge, that is, according to the event
    \[
        A_r^{ij}:=\{\theta_r(i)=j\}.
    \]
    By symmetry among the labels other than $i$, $\E[(\Delta M_r^i)^2\one_{\theta_r(i)=\ell}]$ is the same for every $\ell\ne i$. Hence,
    \[
    \|\Delta M_r^i\|_{L^2}^2=\sum_{\ell: \ell\neq i}\E[\one_{\theta_r(i)=\ell}(\Delta M_r^i)^2]=(N-1)\E[\one_{A_r^{ij}}(\Delta M_r^i)^2].
    \] 
    The analogous statement holds for $\Delta L_r^j$. One then has
        \[
        \E[\one_{A_r^{ij}}(\Delta M_r^i)^2]=\frac1{N-1}\|\Delta M_r^i\|_{L^2}^2,
        \qquad
        \E[\one_{A_r^{ij}}(\Delta L_r^j)^2]=\frac1{N-1}\|\Delta L_r^j\|_{L^2}^2.
    \]
   Moreover, $\E[(\Delta M_r^i)^2\mid \theta_r(i)=\ell]=\E[(\Delta M_r^i)^2]$ for any $\ell$ and the similar inequality holds for $\Delta L_r^j$ as well.
    Therefore, by Cauchy--Schwarz,
    \[
        \left|\E[\one_{A_r^{ij}}\Delta M_r^i\,\Delta L_r^j]\right|
        \leq\frac1{N-1}\|\Delta M_r^i\|_{L^2}\|\Delta L_r^j\|_{L^2}.
    \]

    Suppose $N\ge4$. To handle distinct matching edges, define the $\mathcal F_r$-measurable random variables
    \[
        \begin{aligned}
            a_{r,\ell}^i & :=\E[\Delta M_r^i\mid\mathcal F_r,\theta_r(i)=\ell],\qquad \ell\ne i, \\
            b_{r,\ell}^j & :=\E[\Delta L_r^j\mid\mathcal F_r,\theta_r(j)=\ell],\qquad \ell\ne j.
        \end{aligned}
    \]
    Conditionally on $\mathcal F_r$, the partner labels are uniform and the martingale increments have mean zero. Hence
    \[
        \sum_{\ell\ne i}a_{r,\ell}^i=0,\qquad
        \sum_{\ell\ne j}b_{r,\ell}^j=0.
    \]
    Given $\mathcal F_r$ and the respective partners $\ell,m\notin\{i,j\}$ with $\ell\ne m$, the increments $\Delta M_r^i$ and $\Delta L_r^j$ are independent. The event $\{\theta_r(i)=\ell,\theta_r(j)=m\}$ has conditional probability $\frac{1}{(N-1)(N-3)}$ given $\mathcal F_r$. Summing over $\ell,m$ and using the preceding identities yields
    \[
        \begin{aligned}
             & \E[\one_{(A_r^{ij})^c}\Delta M_r^i\,\Delta L_r^j\mid\mathcal F_r]                               \\
             & \quad=\frac1{(N-1)(N-3)}
            \sum_{\substack{\ell,m\notin\{i,j\}                                                                \\\ell\ne m}}a_{r,\ell}^i b_{r,m}^j\\
             & \quad=\frac1{(N-1)(N-3)}\left[
            \left(\sum_{\ell\notin\{i,j\}}a_{r,\ell}^i\right)
            \left(\sum_{m\notin\{i,j\}}b_{r,m}^j\right)
            -\sum_{\ell\notin\{i,j\}}a_{r,\ell}^i b_{r,\ell}^j\right]                                          \\
             & \quad=\frac{a_{r,j}^i b_{r,i}^j-\sum_{\ell\notin\{i,j\}}a_{r,\ell}^i b_{r,\ell}^j}{(N-1)(N-3)}.
        \end{aligned}
    \]
    Conditional Jensen gives
    \[
        \begin{aligned}
            (a_{r,\ell}^i)^2 & \leq\E[(\Delta M_r^i)^2\mid\mathcal F_r,\theta_r(i)=\ell],\qquad \ell\ne i, \\
            (b_{r,\ell}^j)^2 & \leq\E[(\Delta L_r^j)^2\mid\mathcal F_r,\theta_r(j)=\ell],\qquad \ell\ne j.
        \end{aligned}
    \]
    Taking expectations and using the independence of $\theta_r$ and $\mathcal F_r$ together with the preceding symmetry argument, we obtain
    \[
        \begin{aligned}
            \E[(a_{r,j}^i)^2]                    & \leq\|\Delta M_r^i\|_{L^2}^2,
                                                 & \E[(b_{r,i}^j)^2]                    & \leq\|\Delta L_r^j\|_{L^2}^2,       \\
            \sum_{\ell\ne i}\E[(a_{r,\ell}^i)^2] & \leq (N-1)\|\Delta M_r^i\|_{L^2}^2,
                                                 & \sum_{\ell\ne j}\E[(b_{r,\ell}^j)^2] & \leq (N-1)\|\Delta L_r^j\|_{L^2}^2.
        \end{aligned}
    \]
    The second line follows by summing over the $N-1$ possible partner labels.
    Consequently, Cauchy--Schwarz yields
    \[
        \begin{aligned}
            \E|a_{r,j}^i b_{r,i}^j|
             & \leq\|\Delta M_r^i\|_{L^2}\|\Delta L_r^j\|_{L^2},       \\
            \E\left|\sum_{\ell\notin\{i,j\}}a_{r,\ell}^i b_{r,\ell}^j\right|
             & \leq (N-1)\|\Delta M_r^i\|_{L^2}\|\Delta L_r^j\|_{L^2}.
        \end{aligned}
    \]
    Combining the matched and distinct-edge contributions gives
    \[
        \begin{aligned}
            \left|\E[\Delta M_r^i\,\Delta L_r^j]\right|
             & \leq\left(\frac1{N-1}+\frac{N}{(N-1)(N-3)}\right)
            \|\Delta M_r^i\|_{L^2}\|\Delta L_r^j\|_{L^2}                    \\
             & \leq\frac5{N-1}\|\Delta M_r^i\|_{L^2}\|\Delta L_r^j\|_{L^2}.
        \end{aligned}
    \]
    For $N=2$, only the matched contribution occurs, so the same bound holds.

    Finally, the martingale variance identity gives
    \[
        \sum_{r=0}^{k-1}\|\Delta M_r^i\|_{L^2}^2
        =\operatorname{Var}_{f_k}(\phi)-\operatorname{Var}(M_0^i)
        \leq\operatorname{Var}_{f_k}(\phi),
    \]
    and similarly for $\Delta L_r^j$. Combining \cref{eq:shadow-covariance-increments,eq:shadow-increment-covariance} with Cauchy--Schwarz in $r$ and the preceding variance bounds proves \cref{eq:shadow-covariance}.
    The case $k=0$ follows directly from independent initialization.
\end{proof}

Introduce the shadow empirical and leave-one-out measures
\begin{equation}
    \nu_k^N:=\frac1N\sum_{i=1}^N\delta_{W_k^i},
    \qquad
    \nu_k^{N,-i}:=\frac1{N-1}\sum_{\ell\neq i}\delta_{W_k^\ell}.
    \label{eq:shadow-empiricals}
\end{equation}

\begin{corollary}\label{cor:shadow-empirical}
    Under the assumptions of \Cref{thm:main}, there exists \(C_{T,d,q,\Lambda}\), independent of \(N\) and \(\tau\), such that for every \(i=1,\ldots,N\),
    \begin{equation}\label{eq:shawdowempiricals}
        \begin{aligned}
            \sup_{k\tau\leq T}\E\Wtwo^2(\nu_k^N,f_k)
             & \leq C_{T,d,q,\Lambda}M_q(f_0)^{2/q}\varepsilon_d(N), \\
            \sup_{k\tau\leq T}\E\Wtwo^2(\nu_k^{N,-i},f_k)
             & \leq C_{T,d,q,\Lambda}M_q(f_0)^{2/q}\varepsilon_d(N).
        \end{aligned}
    \end{equation}
\end{corollary}

\begin{proof}
    Fix $k$ with $k\tau\leq T$ and let $\nu$ denote either $\nu_k^N$ or $\nu_k^{N,-i}$, with $\mu=f_k$. Write $\nu=n^{-1}\sum_{\ell\in I}\delta_{W_k^\ell}$, where $I$ is the corresponding set of labels and $n=|I|\in\{N,N-1\}$. Since each particle has law $\mu$, for every Borel set $B$ we have $\E\nu(B)=\mu(B)$ and $\operatorname{Var}(\one_B(W_k^\ell))=\mu(B)(1-\mu(B))$. Applying \cref{prop:shadow-covariance} with $\phi=\psi=\one_B$ gives
    \[
        \left|\operatorname{Cov}\bigl(\one_B(W_k^\ell),\one_B(W_k^{\ell'})\bigr)\right|
        \leq \frac{5}{N-1}\mu(B)(1-\mu(B)),
        \qquad \ell\ne\ell'.
    \]
    Expanding the variance of the empirical average, we therefore obtain
    \begin{equation}\label{eq:shadow-cell-variance}
        \begin{aligned}
            &\E|\nu(B)-\mu(B)|^2\\
             & \quad=\frac1{n^2}\sum_{\ell\in I}
            \operatorname{Var}\bigl(\one_B(W_k^\ell)\bigr)+\frac1{n^2}\sum_{\substack{\ell,\ell'\in I              \\\ell\ne\ell'}}
            \operatorname{Cov}\bigl(\one_B(W_k^\ell),\one_B(W_k^{\ell'})\bigr) \\
             & \quad\leq \frac{\mu(B)(1-\mu(B))}{n^2}
            \left(n+\frac{5n(n-1)}{N-1}\right)                                                   \\
             &\quad \leq \frac{6\mu(B)}n
            \leq \frac{C\mu(B)}N.
        \end{aligned}
    \end{equation}
    Here we used $(n-1)/(N-1)\leq1$ and $n\geq N/2$. The bound is uniform over the omitted label.

    The case $M_q(f_0)=0$ is immediate. Otherwise, apply the common rescaling $v\mapsto v/M_q(f_0)^{1/q}$ to $\mu$ and $\nu$, retaining their notation. The estimate \cref{eq:shadow-cell-variance} remains valid under this rescaling. By \cref{lem:moment}, $M_q(\mu)\leq C_{T,d,q,\Lambda}$ uniformly for $k\tau\leq T$. Set
    \[
        S_0=(-1,1]^d,\qquad
        S_m=(-2^m,2^m]^d\setminus(-2^{m-1},2^{m-1}]^d,\quad m\ge1,
    \]
    and write $b_m=\mu(S_m)$. The moment bound implies $b_m\leq\min\{1,C_{T,d,q,\Lambda}2^{-mq}\}$. At refinement level $j\ge0$, intersect $S_m$ with the dyadic partition of $(-2^m,2^m]^d$ into cubes of side length $2^{m+1-j}$, and denote the resulting partition of $S_m$ by $\mathcal P_{m,j}$. It has at most $2^{dj}$ cells, with $\mathcal P_{m,0}=\{S_m\}$.

    For each $B\in\mathcal P_{m,j}$, Cauchy--Schwarz and \cref{eq:shadow-cell-variance} give
    \[
        \E|\nu(B)-\mu(B)|
        \leq\bigl(\E|\nu(B)-\mu(B)|^2\bigr)^{1/2}
        \leq\frac C{\sqrt N}\sqrt{\mu(B)}.
    \]
    Since $\sum_{B\in\mathcal P_{m,j}}\mu(B)=b_m$ and $\#\mathcal P_{m,j}\leq2^{dj}$, summing over the cells and applying Cauchy--Schwarz again yields
    \[
        \begin{aligned}
            \E\sum_{B\in\mathcal P_{m,j}}|\nu(B)-\mu(B)|
             & \leq\frac C{\sqrt N}\sum_{B\in\mathcal P_{m,j}}\sqrt{\mu(B)}                                      \\
             & \leq\frac C{\sqrt N}\sqrt{\#\mathcal P_{m,j}}\left(\sum_{B\in\mathcal P_{m,j}}\mu(B)\right)^{1/2} \\
             & \leq C2^{dj/2}\sqrt{\frac{b_m}N}.
        \end{aligned}
    \]
    Since $\nu(B),\mu(B)\geq0$ and the cells partition $S_m$, we also have
    \[
        \begin{aligned}
            \E\sum_{B\in\mathcal P_{m,j}}|\nu(B)-\mu(B)| \leq\E\sum_{B\in\mathcal P_{m,j}}\bigl(\nu(B)+\mu(B)\bigr) =\E\nu(S_m)+\mu(S_m)=2b_m,
        \end{aligned}
    \]
    where the last equality uses $\E\nu(S_m)=\mu(S_m)=b_m$. Applying \cite[Lemmas~5--6]{FournierGuillin2015} with $p=2$ gives, for each realization of $\nu$,
    \[
        \Wtwo^2(\nu,\mu)
        \leq C_d\sum_{m\ge0}2^{2m}\sum_{j\ge0}2^{-2j}
        \sum_{B\in\mathcal P_{m,j}}|\nu(B)-\mu(B)|.
    \]
    Taking expectations and applying Tonelli's theorem, we obtain from the preceding two bounds
    \begin{equation}\label{eq:shadow-dyadic-transport}
        \begin{aligned}
            \E\Wtwo^2(\nu,\mu)
             & \leq C_d\sum_{m\ge0}2^{2m}\sum_{j\ge0}2^{-2j}\E\sum_{B\in\mathcal P_{m,j}}|\nu(B)-\mu(B)|       \\
             & \leq C_d\sum_{m\ge0}2^{2m}\sum_{j\ge0}2^{-2j}\min\left\{b_m,2^{dj/2}\sqrt{\frac{b_m}N}\right\}.
        \end{aligned}
    \end{equation}
    Here numerical constants have been absorbed into $C_d$. The term $j=0$ includes the shell mass discrepancy.

    For $0\leq b\leq1$, summation over $j$ gives
    \[
        \sum_{j\ge0}2^{-2j}\min\left\{b,2^{dj/2}\sqrt{\frac bN}\right\}
        \leq C_d
        \begin{cases}
            \sqrt{b/N},             & 2\leq d<4, \\
            \sqrt{b/N}\,[1+\log N], & d=4,       \\
            b^{1-2/d}N^{-2/d},      & d>4.
        \end{cases}
    \]
    Indeed, for $d<4$ the geometric series obtained from the second bound converges. For $d\ge4$ and $Nb\ge1$, split the series at $2^{dj}\asymp Nb$. If $Nb<1$, the bound $C_db$ is no larger than the corresponding expression above.

    The sum over $m$ in \cref{eq:shadow-dyadic-transport} converges because $q>4$ implies $q/2>2$ and, for $d>4$, $q(1-2/d)>2$. Thus $\E\Wtwo^2(\nu,\mu)\leq C_{T,d,q,\Lambda}\varepsilon_d(N)$ for the rescaled measures, uniformly for $k\tau\leq T$ and over the omitted label. Restoring the velocity scale multiplies the squared transport cost by $M_q(f_0)^{2/q}$ and proves both estimates in \cref{eq:shawdowempiricals}.
\end{proof}

\subsection{Proof of the large particle limit}

We now finish the proof of \Cref{thm:main}. In particular, we derive a recursion for the mean-square coupling error between the particle system $\bm V_k$ and the shadow system $\bm W_k$ and then close the estimate using the properties of the shadow system above.

\begin{proof}[Proof of \Cref{thm:main} ]
    Set
    \begin{equation}
        D_k^i:=V_k^i-W_k^i,\qquad E_k^i:=V_k^{\theta_k(i)}-W_k^{i,*}. \label{eq:D-E}
    \end{equation}
    Define
    \begin{align}
        e_k   & :=\frac1N\sum_{i=1}^N\E|D_k^i|^2=\frac1N\sum_{i=1}^N\E|V_k^i-W_k^i|^2,
        \label{eq:e-k}                                                                 \\
        \xi_k & :=\frac1N\sum_{i=1}^N \E\Wtwo^2(f_k^{N,-i}, f_k).
        \label{eq:xi-k}
    \end{align}

    Recalling \Cref{eq:finite-update,eq:original_shadow}, one has
    \begin{equation}
        D_{k+1}^i=A_{P_k^i}D_k^i+B_{P_k^i}E_k^i.\label{eq:D-update}
    \end{equation}
    By \Cref{eq:quadratic-identity}, one has the conditional identity
    \begin{equation}
        \E\left[|D_{k+1}^i|^2\mid\mathcal F_k\vee \sigma( \theta_k, U_k^i)\right]=a_\tau|D_k^i|^2+b_\tau|E_k^i|^2,\label{eq:individual-error}
    \end{equation}
    so that one has
    \begin{equation}
        e_{k+1}=a_\tau e_k+b_\tau\xi_k. \label{eq:exact-e-recursion}
    \end{equation}

    Recalling the shadow empiricals in \Cref{eq:shadow-empiricals}, one has
    \[
        \frac{1}{N}\sum_i \E W_2^2(f_k^{N,-i},\nu_k^{N,-i})
        \leq \frac{1}{N}\sum_i \frac1{N-1}\sum_{\ell\neq i}\E|D_k^\ell|^2=e_k.
    \]
    Then, one has by triangle inequality and \Cref{eq:shawdowempiricals} that
    \begin{equation*}
        \begin{split}
            \xi_k & \le \frac{1}{N}\sum_{i=1}^N 2(\E W_2^2(f_k^{N,-i},\nu_k^{N,-i})
            +\E\Wtwo^2(\nu_k^{N,-i}, f_k))                                          \\
                  & \le 2e_k+2C_{T,d,q,\Lambda}M_q(f_0)^{2/q}\varepsilon_d(N)
        \end{split}
    \end{equation*}
    \phantomsection\label{proof:coupling-recursion}%
    Recall that $a_\tau+b_\tau=1$
    and $b_\tau=O(\tau)$,  one then has for $0\leq k<K$,
    \begin{equation*}
        e_{k+1}\leq(1+b_\tau)e_k+2b_\tau C_{T,d,q,\Lambda}M_q(f_0)^{2/q}\varepsilon_d(N).
    \end{equation*}
    Since $e_0=0$ and recalling \Cref{eq:b-small}, one has
    \begin{align}
        e_k\leq2C_{T,d,q,\Lambda}M_q(f_0)^{2/q}\bigl(e^{2\Lambda(d-1)T}-1\bigr)\varepsilon_d(N).\label{eq:e-rate}
    \end{align}

    Finally, the labelwise coupling of $f_k^N$ and $\nu_k^N$, together with the squared triangle inequality, gives
    \begin{align}\label{eq:final-triangle}
        \E\Wtwo^2(f_k^N, f_k) & \leq2\E\Wtwo^2(f_k^N,\nu_k^N)+2\E\Wtwo^2(\nu_k^N, f_k)\notag \\
                              & \leq2e_k+2\E\Wtwo^2(\nu_k^N, f_k).
    \end{align}
    This together with \Cref{eq:shawdowempiricals} again then gives the final estimate
    \[
        \max_{0\leq k\leq K}\E\Wtwo^2(f_k^N,f_k)\leq C\,M_q(f_0)^{2/q}\varepsilon_d(N).
    \]
    This proves \Cref{eq:main-uniform-bound}, with $C=C(T,d,q,\Lambda)$ independent of $N$, $K$, and $\tau=T/K$.

\end{proof}

\begin{corollary}[Three-dimensional case]\label[corollary]{cor:particle-limit-d3}
    Under the assumptions of \Cref{thm:main}, if $d=3$, then
    \begin{equation}\label{eq:main-uniform-bound-d3}
        \max_{0\le k\le K}\E\Wtwo^2(f_k^N,f_k)
        \le C_{T,q,\Lambda}M_q(f_0)^{2/q}N^{-1/2},
    \end{equation}
    with a constant independent of $N$, $K$, and $\tau$.
\end{corollary}
\begin{proof}
    Apply \Cref{thm:main} with $\varepsilon_3(N)=N^{-1/2}$.
\end{proof}

\section{Entropic convergence of the limit scheme to the Landau equation}\label{sec:entropic}

In this section, we prove \Cref{thm:step2-entropy,cor:errorinW} for Maxwell molecules ($\gamma=0$) under \Cref{ass:step2-initial}. 
Let $\widetilde f$ denote the marginal density given in \Cref{alg:sbm-mean-field}, or equivalently
for $t\in [t_k, t_{k+1}]$ and $Y=(v,v_*)$,
\begin{gather}\label{eq:step2-interpolation}
    \widetilde{f}(t)=\int F(t, Y)\dd v _* =\int e^{(t-t_k) \cL^*}[f_k \otimes f_k] \dd v _*.
\end{gather}
Clearly, $\tilde{f}(t)$ is continuous in time $t$ and $\tilde{f}(t_k)=f_k$.

Define the relative Fisher information
\begin{equation}
    \cI(\tilde{f}|f)=\int \tilde{f} \left|\nabla\log \frac{\tilde{f}}{f}\right|^2\,\dd v,
\end{equation}
and the initial mean and covariance by
\[
    v_c:=\int_{\R^d}v f_0(v)\,\dd v,\qquad
    \Sigma_0:=\int_{\R^d}(v-v_c)\otimes(v-v_c)f_0(v)\,\dd v.
\]
Set
\begin{equation}\label{eq:step2-kappa0}
    \kappa_0:=\Lambda\bigl(\tr\Sigma_0-\lambda_{\max}(\Sigma_0)\bigr).
\end{equation}
As we shall see, under  \Cref{ass:step2-initial}, $\kappa_0$ is strictly positive. Below, we will prove the following more detailed version of \Cref{thm:step2-entropy}, together with \Cref{cor:errorinW}.
\begin{theorem}\label{thm:entropy2}
    Under \Cref{ass:step2-initial}, let $f(t,v)$ solve \Cref{eq:landau}
    with $\gamma=0$ and initial density $f_0$.  The interpolation
    \Cref{eq:step2-interpolation}, initialized with the same $f_0$, satisfies
    \begin{equation}\label{eq:step2-continuous-entropy}
        \sup_{0\le t\le T}\left[\cH(\tilde{f} \mid f)
            +\frac{\kappa_0}{2}\int_0^t\cI(\tilde{f}(t',\cdot)\mid f(t',\cdot))\dd t'\right]
        \le C_T\tau^2,
    \end{equation}
    where $\kappa_0$ is defined in \Cref{eq:step2-kappa0} and is positive by \Cref{lem:eigeval}.   The constant is independent of the partition and depends only on $T$, $d$,
    $\Lambda$, and the initial quantities in \Cref{ass:step2-initial}.
\end{theorem}

\subsection{The proof of the main results}

For $t\in(t_k,t_{k+1})$, \Cref{eq:step2-interpolation} gives
\begin{gather}
    \partial_t \tilde{f}=\int (\nabla_v-\nabla_{v_*})\cdot(A(v-v_*)(\nabla_v-\nabla_{v_*})F)\,\dd v _*.
\end{gather}
Below we will denote $f:=f(t, v), f_*:=f(t, v_*), \tilde{f}=\tilde{f}(t, v), \tilde{f}_*:=\tilde{f}(t, v_*)$.

We now consider the relative entropy defined in \Cref{eq:relativeentropy} between $\tilde{f}$ and $f$:
\begin{gather}
    H(t):=\cH(\widetilde f|f)
    =\int \widetilde f\log\frac{\widetilde f}{f}\,\dd v.
\end{gather}

Fixing $k$,  taking derivative of the relative entropy for $t\in (t_k, t_{k+1})$ gives (see \Cref{app:entropy-justification} for the justification)
\begin{gather}\label{eq:entropyderivative}
    \frac{\dd}{\dd t}H(t)=-D+R_1+R_2,
\end{gather}
where  the dissipation is
\[
    \begin{split}
        D(t)=D(\tilde{f}|f) & :=\iint \nabla_v \log \frac{\tilde{f}}{f} A(v-v_*)
        \left[(\nabla_v\log \frac{\tilde{f}}{f}-\nabla_{v_*}\log\frac{\tilde{f}_*}{f_*}) \tilde{f}\tilde{f}_*  \right]\,\dd v \dd v _* \\
                            & = \frac12 \iint  (q(v)-q(v_*))\cdot A(v-v_*)(q(v)-q(v_*))\widetilde f \widetilde f_* \dd v  \dd v _*,
    \end{split}
\]
where we have introduced
\[
    q(v):=\nabla_v \log \frac{\tilde{f}}{f}=\nabla_v\log \tilde{f}-\nabla_v\log f.
\]

The first remainder term  $R_1$ is the consistency error produced by the pair correlation defect:
\[
    \begin{split}
        R_1 & = -\iint
        q(v)\cdot A(v-v_*)(\nabla_v-\nabla_{v_*})
        \left( F(v,v_*)-\widetilde f \widetilde f_* \right)
        \dd v \dd v _*                                                              \\
            & = -\frac12 \iint (q(v)-q(v_*)) \cdot A(v-v_*) (\nabla_v-\nabla_{v_*})
        \left(F-\widetilde f\widetilde f_* \right) \dd v  \dd v _* .
    \end{split}
\]
The remaining term $R_2$ is
\[
    R_2=\iint \widetilde f(v)(f(v_*)-\widetilde f(v_*))q(v)\cdot  A(v-v_*) (\nabla_v \log f-\nabla_{v_*}\log f_*)\dd v \dd v _* .
\]

Next, we gather some necessary results to control these terms. We first record a simple fact about the covariance matrix of $\tilde{f}$.
\begin{lemma}\label{lem:covariance}
    Let $v_c$ be the mean of $f_k$ and  $\Sigma$ be the covariance matrix of $f_k$. Setting $s=t-t_k$, then one has
    for $s\in [0, \tau]$ that
    \begin{gather*}
        \Sigma_s:=\Sigma(\tilde{f}_t)=\frac{1}{2}\Sigma+\frac{1}{2d}I\tr\Sigma
        +(\frac{1}{2}\Sigma-\frac{1}{2d}I\tr \Sigma) e^{-8\Lambda d s}.
    \end{gather*}
\end{lemma}
\begin{proof}
    Let $(v, v_*)\sim f_k\otimes f_k$. Set $\bar{v}=(v+v_*)/2$ and $\omega=v-v_*=r e$.
    Then, during the evolution $\bar{v}$ is unchanged and $e_s:=\mathcal{R}_{8\Lambda s}e$ satisfies
    \[
        \mathrm d  e_s
        =\sqrt{8\Lambda}(I_d- e_s e_s^{\mathsf T})
        \,\mathrm d\tilde{B}_s
        -4\Lambda(d-1)e_s\,\mathrm ds,
    \]
    Then,
    \[
        \frac{\dd}{\dd s}\E(e_s\otimes e_s)=-8\Lambda d\E(e_s\otimes e_s)+8\Lambda I.
    \]
    Hence,
    \[
        \E(e_s\otimes e_s)=e_0\otimes e_0\exp(-8\Lambda d s)+\frac{1}{d}I(1-e^{-8d\Lambda s}).
    \]
    Hence, the covariance matrix is
    \[
        \Sigma_s=\E(\bar{v}-v_c+\frac{r}{2}e_s)^{\otimes 2}
        =\E (\bar{v}-v_c)^{\otimes 2}+\E \frac{r^2}{4}(e_0\otimes e_0\exp(-8\Lambda d s)+\frac{1}{d}I(1-e^{-8d\Lambda s}))
    \]
    Moreover,
    \[
        \E r^2 e_0\otimes e_0=\E (v-v_*)\otimes (v-v_*)=2\Sigma,
        \quad \E r^2=\E |v-v_*|^2=2M_2-2|v_c|^2.
    \]
    Hence, one has
    \begin{gather}\label{eq:paircovariance}
        \Sigma_s=\frac{1}{2}\Sigma+\frac{1}{2d}I(M_2-|v_c|^2)
        +(\frac{1}{2}\Sigma-\frac{1}{2d}I(M_2-|v_c|^2)) e^{-8\Lambda d s}.
    \end{gather}
\end{proof}

Next, we note the nondegeneracy of the covariance matrix as a consequence of the assumptions.
\begin{lemma}\label{lem:eigeval}
    There is $\kappa >0$ such that every eigenvalue of $\Sigma_0:=\int (v-v_c)\otimes (v-v_c)f_0\,\dd v $ satisfies
    \[
        \lambda_i\ge 1/\cI_0(f_0) \ge  \kappa,\quad  \cI_0(f):=\int |\nabla \log f|^2 f
    \]
    Moreover, the constant $\kappa_0$ in \Cref{eq:step2-kappa0} satisfies
    \begin{gather}
        \kappa_0=\Lambda\bigl(\tr \Sigma_0-\lambda_{\max}(\Sigma_0)\bigr)\ge (d-1)\kappa \Lambda.
    \end{gather}
\end{lemma}

\begin{proof}
    For every unit vector $\xi$, integration by parts gives
    \[
        \int(\xi\cdot (v-v_c))(\xi\cdot\nabla\log f_0(v))f_0(v)\dd v=-1.
    \]
    Taking the square and applying Cauchy--Schwarz, one has
    \[
        \begin{aligned}
            1 & \le (\xi^T\Sigma_0\xi)
            \int|\xi\cdot\nabla\log f_0|^2f_0
            \le (\xi^T\Sigma_0\xi)\cI_0(f_0)                                       \\
              & \le (\xi^T\Sigma_0\xi)K_0^2\int\langle v\rangle^{2m_0}f_0(v)\dd v,
        \end{aligned}
    \]
    where the calculation can be justified easily by the given assumption.
    Thus every eigenvalue satisfies
    \begin{equation}\label{eq:kappa-data}
        \lambda_i\ge \frac{1}{\cI_0(f_0)}\ge \frac{1}{K_0^2\int\langle v\rangle^{2m_0}f_0(v)\dd v}=:\kappa.
    \end{equation}
    The estimate for $\kappa_0$ then immediately follows.
\end{proof}

In \Cref{app:step2-reference}, we will establish the following bounds for logarithmic density.
\begin{lemma}\label{lem:referencesol}
    The solution $f$ to the Landau equation satisfies
    the moment control for $2\le p\le p_*$
    \begin{gather}
        \sup_{t\le T}M_p(f)\le C_{T} M_p(f_0),
    \end{gather}
    and the logarithmic regularity estimate with $m=m_0(m_0+1)$
    \begin{gather}
        \sup_{t\le T} (|\nabla\log f|^2+|\nabla^2 \log f|)\le C_T\langle v\rangle^m
    \end{gather}
\end{lemma}

Define the quantity for a density $\rho$
\begin{equation}\label{eq:Nq}
    \mathcal{N}_q(\rho)=1+\sum_{j=1}^3\int\rho|D_j(\rho)|^{q/j},\quad
    D_j(\rho)=\frac{\nabla^j\rho}{\rho}.
\end{equation}
Later, in \Cref{subsec:localtruncation}, we will establish the following estimate for the solution of the scheme.
\begin{lemma}\label{lem:step2-propagation}
    Set $s=t-t_k$.   For $q>6$, the one-step interpolation \(\tilde f_s\) of Algorithm 3 in the Maxwell case starting from a positive probability density $f_k$ with $\mathcal N_q(f_k)<\infty$, where derivatives are understood weakly, satisfies
    for $s\in [0, \tau]$ that
    \begin{equation}\label{eq:step2-one-propagation}
        \mathcal N_q(\tilde{f}_s)\le(1+C_{d, q,\Lambda}s)\mathcal N_q(f_k),\qquad s\ge0.
    \end{equation}
    Consequently,
    \begin{equation}\label{eq:step2-propagation}
        \sup_{k: k\tau\le T}\mathcal N_q(f_k)
        \le e^{C_{d, q,\Lambda}T}\mathcal N_q(f_0).
    \end{equation}
\end{lemma}

We will need the following coercivity estimate. We attach a short proof in \Cref{app:coercivity}.
\begin{lemma}[Relative Fisher coercivity]\label{lem:relativecoercivity}
    Under the hypotheses of \Cref{ass:step2-initial}, for the constant $\kappa_0$, one has
    \begin{equation}\label{eq:coercivity}
        \quad D(\tilde{f}\mid f)\ge \kappa_0~\cI(\tilde{f}\mid f)-C_0\cH(\tilde{f}\mid f).
    \end{equation}
    Here $C_0$ depends on $T$ and uniform $(2m+4)$-moments of $\tilde{f}$ and $f$.
\end{lemma}

The $R_2$ term can be controlled by
\begin{lemma}\label{lem:prooflem2}
    Suppose \Cref{ass:step2-initial} holds, then for some constant $C_T>0$ independent of $\tau$,
    it holds that
    \[
        |R_2(t)|\le C_T H(t), \qquad 0\le t\le T.
    \]
\end{lemma}

\begin{proof}
    We first consider $R_2$. Noting the identity
    \[
        \tilde{f}q=\nabla(\tilde{f}-f)-(\tilde{f}-f)\nabla\log f,
    \]
    one can then perform integration by parts to deduce
    \begin{gather*}
        R_2=\iint (\tilde{f}_*-f_*)(\tilde{f}-f) K^{R_2}_f(v, v_*)\,\dd v \dd v _*,
    \end{gather*}
    with
    \begin{multline*}
        K^{R_2}_f(v, v_*)=\nabla_v\cdot(A(v-v_*) (\nabla \log f-\nabla_{v_*}\log f_*))+\\
        \nabla\log f\cdot A(v-v_*)(\nabla\log f-\nabla\log f_*).
    \end{multline*}
    One may symmetrize $K^{R_2}_f(v, v_*)$. Direct bound gives
    \[
        |K^{R_2}_f(v, v_*)|\le Cw(v)w(v_*),\quad w(v)=(1+|v|^2)\bigl(1+|\nabla\log f|^2+|\nabla^2\log f|\bigr).
    \]
    Hence, one has (see also \cite{bolley2005weighted})
    \begin{equation*}
        \begin{split}
            |R_2| & \le C\left(\int w(v)|\tilde f-f|\dd v\right)^2                    \\
                  & \le 2C\int w^2(v)(\tilde f+f)\dd v
            \int\bigl(\sqrt{\tilde f}-\sqrt f\bigr)^2\dd v                            \\
                  & \le 2C\left(\int w^2(v)(\tilde f+f)\dd v\right)H(\tilde f\mid f).
        \end{split}
    \end{equation*}
    By \Cref{lem:referencesol} and the moment estimate in \Cref{lem:moment}, also applied at intermediate times, one has $w(v)\le C_T\langle v\rangle^{m+2}$ and a uniform bound on the moments of order $2m+4\le p_*$. Hence, one has the desired bound.
\end{proof}

Define
\begin{equation}\label{eq:flux}
    E_{s}(v)=\int A(v-v_*)\cD(F_s-\tilde{f}\tilde{f}_*)\dd v_*,\qquad \cD=\nabla_v-\nabla_{v_*}.
\end{equation}
and
\[
    \Theta_{f_k}:=\sup_{s\in (0, \tau]}\frac{1}{s^2}\int \frac{1}{\tilde{f}}|E_{s}(v)|^2\,\dd v .
\]

Direct H\"older's inequality in $R_1$ yields
\begin{equation*}
    \begin{split}
        |R_1| & \le \sqrt{\mathcal I(\tilde{f}|f)} \left(\int \frac{1}{\tilde{f}}|E_{s}(v)|^2\,\dd v \right)^{1/2} \\
              & \le \delta \mathcal I(\tilde{f}|f)+C(\delta)\int \frac{1}{\tilde{f}}|E_{s}(v)|^2\,\dd v
        \le \delta \mathcal I(\tilde{f}|f)+C(\delta)\Theta_{f_k}\tau^2.
    \end{split}
\end{equation*}

In \Cref{subsec:localtruncation}, we will use \Cref{lem:step2-propagation} to show that
\begin{gather}
    \sup_{k: k\tau\le T}\Theta_{f_k} \le C_T
\end{gather}
for some constant $C_T$ independent of $\tau$. Hence, for every $k=0,\ldots,K-1$, one then has
\begin{gather}\label{eq:R1bound}
    |R_1(t)|\le \delta~\cI(\tilde{f}(t)| f(t))+C_{T,\delta}\tau^2,\qquad t\in(t_k,t_{k+1}).
\end{gather}

With the above results,  we prove the main results.
\begin{proof}[Proof of \Cref{thm:entropy2}]

    By \Cref{lem:referencesol} and \Cref{lem:step2-propagation}, the calculation and integration by parts to obtain
    \Cref{eq:entropyderivative} can be justified.
    This shows that $t\mapsto H(t)$ is continuous
    and piecewise differentiable.

    By \Cref{lem:relativecoercivity}, \Cref{lem:prooflem2} and \Cref{eq:R1bound}, one finds that for some constant $C_T$ that
    \[
        H'(t)\le -\kappa_0 \cI(\tilde{f}| f)+\delta \cI(\tilde{f}|f)+ C(T,\delta)\tau^2+C_TH(t).
    \]
    This inequality holds for every $t\in (t_k, t_{k+1})$, for each $k=0,\ldots,K-1$, and $e^{-C_Tt}H(t)$ is continuous on $[0, T]$. Noting $H(0)=0$, choosing $\delta=\kappa_0/2$ and applying Gr\"onwall's inequality on each time interval, using continuity at the mesh points, one has
    \[
        H(t)+\frac{\kappa_0}{2}\int_0^t \exp(C_T(t-t'))\cI(t')\dd t' \le C(T,\delta)\tau^2 \int_0^{t} \exp(C_T(t-t')) \dd t'.
    \]
    Since $\exp(C(t-t'))\ge 1$, one then taking the supremum over $t\in [0, T]$ gives the result.
\end{proof}

\begin{proof}[Proof of \Cref{cor:errorinW}]
    Fix $t\in[0,T]$ and write $g=\tilde f(t,\cdot)$ and $f=f(t,\cdot)$. Conservation of momentum and energy gives a common mean $v_c$ and a common centered energy
    \[
        \int |v-v_c|^2g(v)\dd v=\int |v-v_c|^2f(v)\dd v=\mathsf e,\qquad \mathsf e:=\tr\Sigma_0.
    \]
    In the Kantorovich--Rubinstein dual representation, normalize each $1$-Lipschitz test function by $\varphi(v_c)=0$, so that $|\varphi(v)|\le |v-v_c|$. Cauchy--Schwarz and $(\sqrt g+\sqrt f)^2\le2(g+f)$ then give
    \begin{align*}
        W_1(g,f)
         & \le\int |v-v_c|\,|g-f|\dd v                                  \\
         & \le\left(\int |v-v_c|^2(\sqrt g+\sqrt f)^2\dd v\right)^{1/2}
        \left(\int(\sqrt g-\sqrt f)^2\dd v\right)^{1/2}                 \\
         & \le 2\sqrt{\mathsf e\,\cH(g\mid f)}.
    \end{align*}
    Here $\int(\sqrt g-\sqrt f)^2\dd v\le\cH(g\mid f)$ follows by integrating $(\sqrt r-1)^2\le r\log r-r+1$ against $f\dd v$, with $r=g/f$, and using $\int g=\int f=1$. Hence \Cref{thm:step2-entropy} yields
    \[
        \sup_{0\le t\le T}W_1^2(\tilde f(t,\cdot),f(t,\cdot))
        \le 4\mathsf e\sup_{0\le t\le T}\cH(\tilde f(t,\cdot)\mid f(t,\cdot))\le C_T\tau^2.
    \]
    Under the uniform log-Sobolev hypothesis, the transport inequality recalled above and \Cref{thm:step2-entropy} also give $W_2^2(f_k,f(t_k,\cdot))\le C_T\tau^2$.
    Since $W_p\le W_2$, the squared triangle inequality and \Cref{thm:main} yield
    \[
        \E W_p^2(f_k^N,f(t_k,\cdot))
        \le 2\E W_2^2(f_k^N,f_k)+2W_p^2(f_k,f(t_k,\cdot))
        \le C_T\bigl(\varepsilon_d(N)+\tau^2\bigr).
    \]
    For $d=3$, use $\varepsilon_3(N)=N^{-1/2}$.
\end{proof}

\subsection{The local  integral consistency}\label{subsec:localtruncation}

Fix a step $k$ and write $s=t-t_k$, $F_s=e^{s\cL^*}(f_k\otimes f_k)$ and $\tilde f_s=\tilde f(t_k+s)$. Set $D_0(h)=1$ and extend $D_j(h)=\nabla^j h/h$ to pair densities. Constants below may depend on the fixed $d\ge2$, $q$ and $\Lambda$.

In this subsection, we show that $F$ is close to $\tilde{f}\tilde{f}_*$. In particular, we control
$\Theta_{f_k}$. For vector fields $h(v)$ set
\begin{gather}
    \|h\|_s^2=\int |h(v)|^2/\tilde{f}_s(v)\dd v,
\end{gather}
and $\|\cdot\|_s$ gives a norm.
Then,
\[
    \Theta_{f_k}=[\sup_{0<s\le \tau} (\|E_s\|_s/s)]^2.
\]
Moreover, since $E_0=0$,
\[
    \|E_s\|_s=\|\int_0^s \partial_u E_u \dd u\|_s\le \int_0^s \|\partial_u E_u\|_s\,\dd u.
\]

To estimate this quantity, we need some properties for the spherical diffusion.  Below, we set $Y=(v, v_*)$ and consider the equation for $F$.
\begin{lemma}\label{lem:evolutionderilog}
    The transformation $T_s: Y\to Y_s$ is almost surely an orthogonal transformation. Moreover, for any positive probability density $H\in W^{3,1}(\R^{2d})$, with derivatives understood weakly, one has $H_s(Y):=(e^{s\cL^*}H)(Y)$, then for any $b>0$, $j\in\{1,2,3\}$ and $s\in [0, \tau]$,
    \[
        \iint (1+|Y|^2)^b \frac{|\nabla_Y^j H_s|^2}{H_s}\dd Y\le \iint (1+|Y|^2)^b \frac{|\nabla_Y^j H|^2}{H}\dd Y,
    \]
    where
    \[
        |\nabla_Y^j H|^2:=\sum_{|\alpha|=j}\frac{j!}{\alpha!}|\partial^{\alpha} H|^2.
    \]
    More generally, for $a\ge1$, $\ell\ge0$, and $j=1,2,3$,
    \begin{equation}\label{eq:pair-perspective}
        \iint\langle Y\rangle^\ell H_s|D_j(H_s)|^a\dd Y
        \le\iint\langle Y\rangle^\ell H|D_j(H)|^a\dd Y.
    \end{equation}
\end{lemma}
\begin{proof}
    The verification that $T_s$ is orthogonal is straightforward and we omit. Then, since $T_s$ is orthogonal, one has
    \[
        H_s(Y)=\E[H(T_s^{-1}Y)].
    \]
    Here the expectation is taken over the randomness in $T_s$ and this is independent of $Y$. This formula is specific to the Maxwellian case.

    Since $(a, z)\to |z|^2/a$ is convex, the Jensen's inequality implies that
    \[
        \iint (1+|Y|^2)^b \frac{|\nabla_Y^j H_s|^2}{H_s}\dd Y
        \le \E \iint (1+|Y|^2)^b \frac{|\nabla_Y^j H(T_s^{-1}Y)|^2}{H(T_s^{-1}Y)}\dd Y.
    \]
    Since $T_s$ is orthogonal, one has
    \[
        |\nabla_Y^j H(T_s^{-1}Y)|=|(\nabla^j H)(T_s^{-1}Y)|.
    \]
    Change of variables then shows that
    \[
        \iint (1+|Y|^2)^b \frac{|(\nabla^j H)(T_s^{-1}Y)|^2}{H(T_s^{-1}Y)}\dd Y=\iint (1+|Y|^2)^b \frac{|\nabla_Y^j H|^2}{H}\dd Y.
    \]

    Replacing $|z|^2/a$ by the convex function $(h,z)\mapsto |z|^a/h^{a-1}$ proves \Cref{eq:pair-perspective} by the same argument.
\end{proof}

We next estimate $\|E_s\|_s$.
\begin{lemma}[Estimate in the time-dependent norm]\label{lem:movingnorm}
    Let $p,r>2$ satisfy $1/p+1/r=1/2$. For smooth $\tilde f_s>0$ and $E_s$ with $E_0=0$, suppose that, uniformly for $0\le s\le\tau$,
    \[
        \|\partial_sE_s\|_s\le A_0,\qquad
        \left\|\frac{\partial_s\tilde f_s}{\tilde f_s}\right\|_{L^p(\tilde f_s)}\le B_0,\qquad
        \left\|\frac{E_s}{\tilde f_s}\right\|_{L^r(\tilde f_s)}\le C_0,
    \]
    where $\|u\|_{L^a(\tilde f_s)}=(\int |u|^a\tilde f_s\dd v)^{1/a}$. Then
    \[
        \|E_s\|_s\le(A_0+B_0C_0/2)s,\qquad 0\le s\le\tau.
    \]
\end{lemma}
\begin{proof}
    Set $X(s)=\|E_s\|_s^2$, $z_s=E_s/\tilde f_s$, and $b_s=(\partial_s\tilde f_s)/\tilde f_s$. Formally,
    \[
        X'(s)=2\int\tilde f_s z_s\cdot\frac{\partial_sE_s}{\tilde f_s}\dd v
        -\int\tilde f_s|z_s|^2b_s\dd v.
    \]
    The relation between $p$ and $r$ gives
    \[
        \|z\|_{L^{2p/(p-1)}(\tilde f_s)}^2
        \le\|z\|_{L^2(\tilde f_s)}\|z\|_{L^r(\tilde f_s)}.
    \]
    Cauchy--Schwarz and H\"older's inequality therefore yield
    \[
        |X'(s)|\le 2A_0\sqrt{X(s)}+B_0\|z_s\|_{L^2(\tilde f_s)}\|z_s\|_{L^r(\tilde f_s)}
        \le(2A_0+B_0C_0)\sqrt{X(s)}.
    \]
    To justify the differentiation, take smooth compactly supported cutoffs $0\le\chi_R\uparrow1$ and replace $X$ by
    \[
        X_{R,\varepsilon}(s)=\int\chi_R\frac{|E_s|^2}{\tilde f_s+\varepsilon}\dd v,
        \qquad z_{s,\varepsilon}=\frac{E_s}{\sqrt{\tilde f_s(\tilde f_s+\varepsilon)}}.
    \]
    The two differentiated terms are bounded, respectively, by
    \[
        2A_0\sqrt{X_{R,\varepsilon}}\quad\hbox{and}\quad B_0C_0\sqrt{X_{R,\varepsilon}},
    \]
    using the same interpolation for the finite measure $\chi_R\tilde f_s\dd v$ and $|z_{s,\varepsilon}|\le|z_s|$. Since $X_{R,\varepsilon}(0)=0$, integrate the resulting inequality for $\sqrt{X_{R,\varepsilon}+\eta}$ and let $\eta\downarrow0$. Finally, let $\varepsilon\downarrow0$ and $R\to\infty$ by monotone convergence to obtain the stated bound.
\end{proof}

Introduce
\begin{equation}\label{eq:Pdef}
    \mathcal{B}(\rho)
    :=\sum_{j=1}^3\iint (1+|Y|^2)^{j+1}
    \frac{|\nabla_Y^j(\rho\otimes\rho)|^2}{\rho\otimes\rho}\dd Y.
\end{equation}
We use the same notation $\cB[F]$ for the right-hand expression with a general pair density $F$ in place of $\rho\otimes\rho$.

The main claim here is the following.
\begin{proposition}\label{prop:Theta}
    For $q>6$ and $P_q=8q/(q-6)$,
    \begin{gather}\label{eq:step2-Theta-bound}
        \Theta_{f_k}\le C(1+M_2(f_k))^2{\bigl[\mathcal B(f_k)+1+M_{P_q}(f_k)+\mathcal N_q(f_k)\bigr]}.
    \end{gather}
\end{proposition}

\begin{proof}
    Introduce
    \[
        (\cK H)(v)=\int A(v-v_*)\cD H(v,v_*)\dd v_*.
    \]
    Then, one has
    \[
        \partial_s E_s=\cK(\cL^* F_s)-\partial_s \cK (\tilde{f}\tilde{f}_*).
    \]
    Since for the Maxwellian case $\nabla\cdot A(\omega)=-\Lambda (d-1)\omega$, one has
    \[
        \cL^*F=-2\Lambda(d-1)(v-v_*)\cdot\cD F+A(v-v_*):\cD^2F.
    \]
    Hence, one has the following pointwise bound
    \begin{align}
        |\cL^* F|         & \le C\bigl((1+|Y|^2)|\nabla_Y^2F|+(1+|Y|^2)^{1/2}|\nabla_YF|\bigr),\label{eq:Lbound} \\
        |\nabla_Y\cL^* F| & \le C\bigl((1+|Y|^2) |\nabla_Y^3F|
        +(1+|Y|^2)^{1/2}|\nabla_Y^2F|+|\nabla_YF|\bigr).\label{eq:gradLbound}
    \end{align}
    Then, H\"older inequality, \Cref{lem:evolutionderilog} and the above imply that
    \begin{multline*}
        \|\cK(\cL^* F_s)\|_s^2
        \le \iint \frac{|A(v-v_*)\cD \cL^* F_s|^2}{F_s} \dd Y \\
        \le C\iint\frac{(1+|Y|^2)^4|\nabla_Y^3F_s|^2
            +(1+|Y|^2)^3|\nabla_Y^2F_s|^2+(1+|Y|^2)^2|\nabla_YF_s|^2}{F_s}\dd Y
        \le C\cB(f_k).
    \end{multline*}
    For the second term,
    \[
        \cK(\tilde{f}\tilde{f}_*)
        =(A*\tilde{f})\nabla\tilde{f}+\tilde{f} \Lambda(d-1)(v-v_c)
    \]
    Note
    \[
        A*\tilde{f}=A(v-v_c)+\Lambda\bigl((M_2-|v_c|^2)I-\Sigma_s\bigr)
    \]
    By \Cref{lem:covariance}, one can find that $|\partial_s \Sigma_s|\le C(1+M_2)$ and thus
    \begin{gather*}
        \|\partial_s \cK(\tilde{f}\tilde{f}_*)\|_s^2
        \le C(1+M_2(f_k))^2\int \frac{|\nabla \tilde{f}|^2}{\tilde{f}}
        +\l v\r^2\frac{|\partial_s\tilde{f}|^2}{\tilde{f}}
        +\l v\r^4\frac{|\partial_s\nabla \tilde{f}|^2}{\tilde{f}} \dd v.
    \end{gather*}
    Since
    \[
        \partial_s \tilde{f}=\int \cL^* F_s \dd v_*,
    \]
    one may then apply the convexity of $(a,z)\to |z|^2/a$ to find that
    \[
        \|\partial_s \cK(\tilde{f}\tilde{f}_*)\|_s^2
        \le C(1+M_2(f_k))^2\mathcal{B}[F_s]\le C(1+M_2(f_k))^2\mathcal{B}(f_k).
    \]
    Combining these estimates yields
    \begin{equation}\label{eq:local-source-norm}
        \sup_{0\le s\le\tau}\|\partial_sE_s\|_s
        \le C(1+M_2(f_k))\mathcal B(f_k)^{1/2}.
    \end{equation}
    For the other two bounds, choose
    \begin{equation}\label{eq:moving-exponents}
        p=\frac{4q}{q+2},\qquad r=\frac{4q}{q-2},\qquad \frac1p+\frac1r=\frac12.
    \end{equation}
    The flux identity above also gives
    \[
        E_s=\int\bigl[A(v-v_*)\cD F_s-(A*\tilde f_s)\nabla_vF_s
            -\Lambda(d-1)(v-v_c)F_s\bigr]\dd v_*.
    \]
    For $u(v)=\int U(v,v_*)\dd v_*$, conditional Jensen's inequality gives
    \[
        \left\|\frac{u}{\tilde f_s}\right\|_{L^a(\tilde f_s)}^a
        \le\iint F_s\left|\frac{U}{F_s}\right|^a\dd Y,\qquad a\ge1.
    \]
    Applying this to $\partial_s\tilde f_s=\int\cL^*F_s\dd v_*$ and to $E_s$, respectively, yields
    \begin{align}
        \left\|\frac{\partial_s\tilde f_s}{\tilde f_s}\right\|_{L^p(\tilde f_s)}^p
         & \le C\iint F_s\bigl(\langle Y\rangle^{2p}|D_2(F_s)|^p
        +\langle Y\rangle^p|D_1(F_s)|^p\bigr)\dd Y,\label{eq:local-density-norm} \\
        \left\|\frac{E_s}{\tilde f_s}\right\|_{L^r(\tilde f_s)}^r
         & \le C(1+M_2(f_k))^r\iint F_s\bigl(\langle Y\rangle^{2r}|D_1(F_s)|^r
        +\langle Y\rangle^r\bigr)\dd Y.\label{eq:local-flux-norm}
    \end{align}
    The weighted estimate in \Cref{lem:Bbound} below applies because
    \[
        \frac{2pq}{q-2p}=\frac{2rq}{q-r}=P_q,\qquad
        \frac{pq}{q-p}=r<P_q.
    \]
    The term without derivatives is controlled by orthogonality and the $P_q$-moment. Thus, together with \Cref{eq:local-source-norm}, we have
    \begin{equation}\label{eq:local-three-norms}
        \begin{split}
            \sup_{0\le s\le\tau}\left\|\frac{\partial_s\tilde f_s}{\tilde f_s}\right\|_{L^p(\tilde f_s)}
             & \le C\bigl(1+M_{P_q}(f_k)+\mathcal N_q(f_k)\bigr)^{1/p},             \\
            \sup_{0\le s\le\tau}\left\|\frac{E_s}{\tilde f_s}\right\|_{L^r(\tilde f_s)}
             & \le C(1+M_2(f_k))\bigl(1+M_{P_q}(f_k)+\mathcal N_q(f_k)\bigr)^{1/r}.
        \end{split}
    \end{equation}
    Applying \Cref{lem:movingnorm}, using $1/p+1/r=1/2$ and $(a+b)^2\le2a^2+2b^2$, proves \Cref{eq:step2-Theta-bound} for smooth data.

    For data with the stated regularity, we apply the preceding argument to $f_k^\varepsilon=f_k*\gamma_\varepsilon^{(d)}$, where $\gamma_\varepsilon^{(n)}$ is the centered Gaussian density with covariance $\varepsilon I_n$. Jensen's inequality gives $\mathcal N_q(f_k^\varepsilon)\le\mathcal N_q(f_k)$ and, with polynomial weights, $\mathcal B(f_k^\varepsilon)\le C\mathcal B(f_k)$ for $0<\varepsilon\le1$; the $P_q$-moments converge. Orthogonality gives
    \[
        F_s^\varepsilon=F_s*\gamma_\varepsilon^{(2d)},\qquad
        \tilde f_s^\varepsilon=\tilde f_s*\gamma_\varepsilon^{(d)}.
    \]
    The assumed normalized derivative and moment bounds imply weighted $W^{3,1}$ convergence with weight $\langle Y\rangle^4$. The marginal formulas and convergence of the covariance coefficients then imply convergence of the fluxes in $L^1$. Lower semicontinuity of $\int|E_s|^2/\tilde f_s$ passes the estimate to the limit. This also justifies the differentiations and integrations used above.
\end{proof}

Hence, it suffices to estimate $\mathcal{B}(f_k)$ and $M_{P_q}(f_k)+\mathcal N_q(f_k)$.
For this purpose,  we first prove \Cref{lem:step2-propagation} and then apply it to estimate
$\mathcal{B}(f_k)$.
\begin{proof}[Proof of \Cref{lem:step2-propagation}]
    We first prove the estimate for smooth data.

    Set $H(x,y)=f_k(x)f_k(y)$. Consider the orthogonal map $T_s: Y\to Y_s$. Let $H_s=H\circ T_s^{-1}$.
    Then $\tilde{f}_s(v)=\E\int H_{s}(v, v_*)\dd v_*$.

    Write
    \[
        (v, v_*)=T_s(x , y)=(A_sx+B_sy, B_sx+A_s y),\quad A_s=\frac{1}{2}(I+P_s)
        , B_s=\frac{1}{2}(I-P_s),
    \]
    where
    $P_s=\mathcal R_{8\Lambda s}$
    Define
    \begin{multline*}
        \mathcal{T}_j(P_s; x, y):=\frac{\nabla_v^j H_s(v,  v_*)}{H_s(v, v_*)}|_{(v, v_*)=T_s(x, y)} \\
        = \sum_{a=0}^j\binom ja\operatorname{Sym}
        \left[(A_s^{\otimes a}D_a(f_k)(x))
        \otimes(B_s^{\otimes(j-a)}D_{j-a}(f_k)(y))\right],
    \end{multline*}
    where for an order $j$ tensor $U$,
    \[
        (\operatorname{Sym}(U))_{i_1,\cdots, i_j}=\frac{1}{j!}\sum_{\pi \in S_j} U_{i_{\pi(1)}\cdots i_{\pi(j)}}.
    \]

    For a given $q$, the quantity we are interested in is
    $$\int \tilde{f} \left|\frac{\nabla_v^j \tilde{f}}{\tilde{f}}\right|^{q/j}\dd v .$$
    Note that $H_s(v, v_*) \mu(dP_s) \dd v \dd v_*$ is the joint distribution of $(P_s, v, v_*)$.
    By Jensen's inequality
    \begin{multline*}
        \left|\frac{\nabla_v^j \tilde{f}}{\tilde{f}}\right|^{q/j}
        =\left|\frac{\int \E[\nabla_v^j H_s(v, v_*)] \dd v _* }{\tilde{f}}\right|^{q/j}
        =\left| \E_{P_s,V_*} \left[\frac{\nabla_v^j H_s(v, V_*)}{H_s(v, V_*)} \,\middle|\, v \right]\right|^{q/j}\\
        \le  \E_{P_s,V_*}\left[\left|\frac{\nabla_v^jH_s(v,V_*)}{H_s(v,V_*)}\right|^{q/j}\,\middle|\,v\right]
        =\frac{1}{\tilde{f}(v)}\E \int H_s(v, v_*) \left| \frac{\nabla_v^j H_s(v, v_*)}{H_s(v, v_*)}\right|^{q/j} \dd v _*
    \end{multline*}
    Integrating this with respect to $\tilde{f} \dd v $ and performing change of variables $(v, v_*)\to (x, y)$, one finds
    \[
        \int \tilde{f} \left|\frac{\nabla_v^j \tilde{f}}{\tilde{f}}\right|^{q/j}\dd v  \le \E
        \iint H(x, y) \left| \mathcal{T}_j(P_s, x, y)\right|^{q/j} \dd x\dd y
    \]
    The randomness is then moved to $\mathcal{T}_j$.  At $s=0$, both sides are equal.

    Defining
    \[
        \Phi_j(P_s)=\iint H(x, y) \left| \mathcal{T}_j(P_s, x, y)\right|^{q/j} \dd x\dd y.
    \]
    Dynkin's formula implies that
    \begin{gather}\label{eq:auxPhievol}
        \begin{split}
            \int \tilde{f} \left|\frac{\nabla_v^j \tilde{f}}{\tilde{f}}\right|^{q/j}\dd v  \le
            \E \Phi_j(P_s) & =\Phi_j(I)+\E \int_0^s \mathcal L_{\mathrm{rot}}\Phi_j(P_u)\,\dd u                                        \\
                           & \le\int f_k \left|\frac{\nabla_v^j f_k}{f_k}\right|^{q/j}\dd v +\|\cL_{\mathrm{rot}} \Phi_j\|_{\infty} s.
        \end{split}
    \end{gather}

    Recall the basic setup in \Cref{subsec:sbmbasics}, the generator of $P_s e$
    for any fixed $e$ is $4\Lambda \Delta_{\mathbb{S}^{d-1}}$. Note that $\cL_{\mathrm{rot}}$ is related to the spherical generator by
    \[
        \mathcal L_{\mathrm{rot}}[R\mapsto u(Re)](R)=4\Lambda(\Delta_{\mathbb S^{d-1}}u)(Re)
    \]
    for any smooth function $u$ on $\mathbb{S}^{d-1}$ and $e\in \mathbb{S}^{d-1}$.
    Hence, one then has
    \[
        \mathcal L_{\mathrm{rot}}=4\Lambda \Delta_{SO(d)}.
    \]
    In particular, recalling $J_{a,b}$ there, one can define a vector field
    \[
        X_{a,b}\Psi(R)=\frac{\dd}{\dd t}\Psi(e^{tJ_{a,b}}R)|_{t=0},\quad
        X_{a,b}=\sum_{c=1}^d (R_{bc}\frac{\partial}{\partial R_{ac}}-R_{ac}\frac{\partial}{\partial R_{bc}}).
    \]
    Then,
    \[
        \cL_{\mathrm{rot}}=4\Lambda \sum_{a<b} X_{a,b}^2.
    \]
    Each vector $X_{a,b}$ will then act on $\Phi_j$ through $A_s^{\otimes a}$ and $B_s^{\otimes (j-a)}$,
    which will result in bounded terms. Hence,
    \[
        |\mathcal T_j|+|X_{ab}\mathcal T_j|+|X_{ab}^2\mathcal T_j|
        \le C_j(\sum_{a=0}^j|D_a(f_k)(x)|\,|D_{j-a}(f_k)(y)|).
    \]
    Recalling that $H(x, y)=f_k(x)f_k(y)$, one has
    \begin{multline*}
        | \cL_{\mathrm{rot}} \Phi_j|\le C_{d,q,\Lambda}\sum_{a<b}\iint f_k(x)f_k(y)\\
        \quad\Bigl( \left| \mathcal{T}_j(P_s, x, y)\right|^{q/j-2}|X_{ab}\mathcal{T}_j(P_s, x, y)|^2
        \\\quad+\left| \mathcal{T}_j(P_s, x, y)\right|^{q/j-1}|X_{ab}^2\mathcal{T}_j(P_s, x, y)|\Bigr)\dd x\dd y\\
        \le C \iint f_k(x)f_k(y)
        (\sum_{a=0}^j|D_a(f_k)(x)|\,|D_{j-a}(f_k)(y)|)^{q/j}\\
        \le C \iint f_k(x)f_k(y)
        \sum_{a=0}^j(|D_a(f_k)(x)|\,|D_{j-a}(f_k)(y)|)^{q/j}
    \end{multline*}
    For $a=0$ or $a=j$, the term in the sum gives the same
    $\int f_k |D_j(f_k)|^{q/j}\dd x$.
    For $1\le a\le j-1$, by Young's inequality $ab\le a^p/p+b^{p'}/p'$, one has
    \begin{multline*}
        \iint f_k(x)f_k(y)  |D_a(f_k)(x)|^{q/j}\,|D_{j-a}(f_k)(y)|^{q/j}\\
        \le C(\iint f_k(x)f_k(y) |D_a(f_k)(x)|^{q/a}
        +\iint f_k(x)f_k(y) |D_{j-a}(f_k)(y)|^{q/(j-a)})
    \end{multline*}
    Hence, one has
    \[
        \|\cL_{\mathrm{rot}} \Phi_j(R)\|_{\infty}
        \le C_{d, q,\Lambda}\sum_{a=1}^j \int f_k(x) |D_a(f_k)(x)|^{q/a}\dd x.
    \]
    Substituting this back into \Cref{eq:auxPhievol} and summing over $j=1,2,3$,
    one then obtains \Cref{eq:step2-one-propagation}.

    For the stated Sobolev data, let $f_k^\varepsilon=f_k*\gamma_\varepsilon^{(d)}$, with $\gamma_\varepsilon^{(n)}$ the centered Gaussian density of covariance $\varepsilon I_n$. Jensen's inequality gives $\mathcal N_q(f_k^\varepsilon)\le\mathcal N_q(f_k)$. Orthogonality and marginalization give $\tilde f_s^\varepsilon=\tilde f_s*\gamma_\varepsilon^{(d)}$; moreover, the rotation representation and Fubini's theorem imply $\tilde f_s\in W^{3,1}$. Hence Gaussian approximation converges in $W^{3,1}$, and lower semicontinuity of the perspective integrals passes \Cref{eq:step2-one-propagation} to the limit. Under \Cref{ass:step2-initial}, H\"older's inequality gives $f_0\in W^{3,1}$, so the estimate applies inductively.

    By setting $s=\tau$, one obtains the recursion from $k$ to $k+1$.
    The last claim then follows directly by this recursion.
\end{proof}

Then, we conclude the control of $\mathcal{B}(f_k)$ in \Cref{eq:Pdef}.
\begin{lemma}\label{lem:Bbound}
    It holds that
    \begin{gather}
        \begin{split}
            \mathcal{B}(f_k) & \le C(\sum_{j=0}^3\int \l v\r^8 f_k |D_j(f_k)|^{2}\dd v)^2 \\
                             & \le C(1+M_{8q/(q-6)}(f_k)+\mathcal{N}_q(f_k))^2.
        \end{split}
    \end{gather}
    In addition, for $j\in\{1,2,3\}$, $a\ge1$, $ja<q$, and $\ell\ge0$ with $\ell q/(q-ja)\le P_q$,
    \begin{equation}\label{eq:pair-weighted-control}
        \sup_{0\le s\le\tau}\iint\langle Y\rangle^\ell F_s|D_j(F_s)|^a\dd Y
        \le C\bigl(1+M_{P_q}(f_k)+\mathcal N_q(f_k)\bigr).
    \end{equation}
\end{lemma}

\begin{proof}
    Direct calculation reveals that
    \[
        \frac{|\nabla_Y^j(\rho\otimes\rho)|^2}{\rho\otimes\rho}
        =\sum_{b=0}^j\binom jb\rho(v)\rho(v_*)
        |D_b(\rho)(v)|^2|D_{j-b}(\rho)(v_*)|^2.
    \]
    Setting $\rho=f_k$ and using the trivial bound
    \[
        (1+|Y|^2)^{j+1}\le \l v\r^8 \l v_*\r^8,
    \]
    one can find that
    \[
        \mathcal{B}(f_k)\le C\left(\sum_{j=0}^3 \int f_k(v)\l v\r^8|D_j(f_k)|^2\right)^2.
    \]
    For $j=1,2,3$, Young's inequality gives
    \[
        \l v\r^8|D_j(f_k)|^2\lesssim \l v\r^{\frac{8q}{q-2j}}+ |D_j(f_k)|^{\frac{q}{j}}.
    \]
    The term $j=0$ is bounded by $C(1+M_8(f_k))$, and $P_q>8$. Since $8q/(q-2j)\le P_q$ for $j=1,2,3$, the result then follows.

    For $H=f_k\otimes f_k$, the product rule and H\"older's inequality in each velocity give
    \[
        \iint H|D_j(H)|^{q/j}\dd Y\le C\mathcal N_q(f_k).
    \]
    Indeed, each mixed product with derivative orders $b$ and $j-b$ is bounded by
    \[
        \left(\int f_k|D_b(f_k)|^{q/b}\right)^{b/j}
        \left(\int f_k|D_{j-b}(f_k)|^{q/(j-b)}\right)^{(j-b)/j};
    \]
    the endpoint terms are immediate. Another application of H\"older gives
    \[
        \begin{split}
            \iint\langle Y\rangle^\ell H|D_j(H)|^a\dd Y
             & \le\left(\iint H|D_j(H)|^{q/j}\dd Y\right)^{ja/q}                             \\
             & \quad\times\left(\iint H\langle Y\rangle^{\ell q/(q-ja)}\dd Y\right)^{1-ja/q} \\
             & \le C\bigl(1+\mathcal N_q(f_k)+M_{P_q}(f_k)\bigr).
        \end{split}
    \]
    Applying \Cref{eq:pair-perspective} proves \Cref{eq:pair-weighted-control}.
\end{proof}

By \Cref{lem:Bbound}, \Cref{lem:step2-propagation}, \Cref{lem:moment,eq:step2-Theta-bound} and the assumptions in \Cref{ass:step2-initial}, one then conclude for some constant $C_T$ that
\begin{equation}\label{eq:step2-uniform-Theta}
    \sup_{k: k\tau\le T}\Theta_{f_k}\le C_T.
\end{equation}

\section*{AI Statement}

The authors used GPT-6 Astra (OpenAI) for mathematical discussions, manuscript drafting, and language editing. The AI tool played important roles in figuring out the covariance estimate in \Cref{prop:shadow-covariance}, and  the detailed local truncation error in \Cref{lem:movingnorm} during the discussion with authors.
The authors take full responsibility for the mathematical arguments, results, and final content of the manuscript.

\section*{Acknowledgements} JAC was supported by the Advanced Grant Nonlocal-CPD (Nonlocal PDEs for Complex Particle Dynamics: Phase Transitions, Patterns and Synchronization) of the European Research Council Executive Agency (ERC) under the European Union Horizon 2020 research and innovation programme (grant agreement No. 883363) and  partially supported by the EPSRC EP/V051121/1. The work of LL is partially supported by the National Natural Science Foundation of China (NSFC) under grant No.~12371400. JAC and LL are partially supported by the Innovation Base for Mathematical Life Sciences and Artificial Intelligence under Number B26039.

\appendix

\section{Uniqueness of the two-particle forward equation}
\label[appendix]{app:pair-forward-uniqueness}

\begin{lemma}[Uniqueness of the two-particle forward equation]
    \label[lemma]{lem:pair-forward-uniqueness}
    Let \(\mathcal L^*\) be the Maxwell operator in \Cref{eq:pair_initial_value_problem}, and fix \(k\) with \(f_k\in\mathcal P_2(\mathbb R^d)\). The initial-value problem
    \[
        \partial_s\zeta_s=\mathcal L^*\zeta_s,
        \qquad \zeta_0=f_k^{\otimes2},
        \qquad 0\leq s\leq \tau,
    \]
    has at most one weak solution in \(C([0,\tau];\mathcal P_2(\mathbb R^{2d}))\), equipped with the \(W_2\) topology.
\end{lemma}

\begin{proof}
    Let \(J_{ab}=e_ae_b^{\mathsf T}-e_be_a^{\mathsf T}\), where \(e_1,\ldots,e_d\) is the standard basis. The process \Cref{eq:candidate_pair} satisfies the Stratonovich system
    \[
        \mathrm dX_s
        =\sqrt{2\Lambda}\sum_{a<b}J_{ab}(X_s-X_s^*)\circ\mathrm d B_s^{ab},
        \qquad
        \mathrm dX_s^*=-\mathrm dX_s,
    \]
    where the \(B^{ab}\) are independent standard Brownian motions. Its diffusion coefficients are linear, and its It\^o drifts are \(-2\Lambda(d-1)(X_s-X_s^*)\) and the opposite vector. Write \(b(z)\) and \(\Sigma(z)\) for the drift vector and diffusion coefficient of the corresponding It\^o equation on \(\mathbb R^{2d}\), where \(z=(v,v_*)\). Both coefficients are linear and globally Lipschitz, so the martingale problem for \(\mathcal L\) is well posed for the initial law \(f_k^{\otimes2}\)\cite[Chapter~32]{Kallenberg2021}. Its diffusion covariance \(a(z):=\Sigma(z)\Sigma(z)^{\mathsf T}\) satisfies
    \[
        |b(z)|+\|a(z)\|_{\mathrm{HS}}\leq C(1+|z|^2),
    \]
    where \(C\) depends only on \(d\) and \(\Lambda\).

    To pass from uniqueness of the martingale problem to uniqueness of the forward equation, let \((\zeta_s)_{0\leq s\leq \tau}\) be any weak solution with initial law \(f_k^{\otimes2}\) in the stated class. Continuity on the compact time interval gives a uniform second-moment bound, and hence
    \[
        \int_0^\tau\int_{\mathbb R^{2d}}\bigl(|b(z)|+\|a(z)\|_{\mathrm{HS}}\bigr)\,\zeta_s(\mathrm dz)\,\mathrm ds
        \leq C\tau\left(1+\sup_{0\leq s\leq \tau}\int_{\mathbb R^{2d}}|z|^2\,\zeta_s(\mathrm dz)\right)<\infty.
    \]
    Since \(W_2\) continuity implies narrow continuity, the superposition principle \cite[Theorem~2.5]{Trevisan2016} applies: there is a probability measure on \(C([0,\tau];\mathbb R^{2d})\) solving the martingale problem for \(\mathcal L\) whose time-\(s\) marginal is \(\zeta_s\). Uniqueness of that martingale problem therefore implies uniqueness of the weak forward equation in the stated class.
\end{proof}

\section{Logarithmic regularity of the reference solution}
\label[appendix]{app:step2-reference}

We record an initial-data criterion for the reference bounds used in the
entropy estimate.  The argument uses the rotation and conditional Gaussian
representation of the Maxwell equation; see \cite{DMN} for this representation.

\begin{lemma}[Reference bounds from the initial density]
    \label{lem:step2-reference}
    Let $f_0$ be a strictly positive $C^2$ probability density on $\R^d$, $d\ge2$, with finite second moment. Suppose, for some $m_0\ge1$ and $K<\infty$,
    \[
        |\nabla_v\log f_0(v)|+|\nabla_v^2\log f_0(v)|
        \le K\langle v\rangle^{m_0}.
    \]
    Set $m=m_0(m_0+1)$.  The Maxwell Landau solution with initial density $f_0$ is positive and smooth for $t>0$, and continuous at $t=0$ in
    $C^2_{\mathrm{loc}}$.  For every finite $T$,
    \begin{equation}\label{eq:step2-reference-bound}
        |\nabla_v\log f(t,v)|^2+|\nabla_v^2\log f(t,v)|
        +\frac{|\nabla_v^2 f(t,v)|}{ f(t,v)}
        \le C_T\langle v\rangle^m,\qquad 0\le t\le T.
    \end{equation}
    For every $p\ge2$ with finite initial $p$th moment,
    \[
        \sup_{t\le T}\int\langle v\rangle^p f(t,v)\dd v
        \le C_{T,p}\int\langle v\rangle^p f_0(v)\dd v.
    \]
    The constants depend only on the stated initial quantities, $T$, $d$ and $\Lambda$.
\end{lemma}

\begin{proof}
    By translation, we may assume $v_c=0$; translating back changes only the constants. Write $\Sigma_t=\int (v-v_c)\otimes (v-v_c) f_t(v)\dd v$, $K=\tr~\Sigma_t=\tr~\Sigma_0$.  Testing the Landau equation against $(v-v_c)\otimes (v-v_c)$ gives
    $\dot\Sigma_t=4\Lambda K I_d-4 d \Lambda\Sigma_t$.  Then,
    \begin{equation}\label{eq:Landaucov}
        \Sigma_t=\frac Kd I_d+e^{-4d\Lambda t}
        \left(\Sigma_0-\frac KdI_d\right).
    \end{equation}
    This should be compared to
    \Cref{eq:paircovariance} where
    \[
        \Sigma(\tilde{f}_s)=\frac{I_d}{d}\tr~\Sigma_0+\frac{1+e^{-8\Lambda d s}}{2} (\Sigma_0-\frac{I_d}{d}\tr~\Sigma_0).
    \]
    It seems that they have different decaying rate. In fact, they match for $s$ small (the error is $O(s^2)$ leading to global convergence). The decaying rate in $\Sigma(\tilde{f}_s)$ is higher but the limit still contains $\Sigma_0/2$ which reflects that the paired center is constant in the two-particle dynamics.

    Every classical
    solution with this initial density has the displayed covariance in \Cref{eq:Landaucov}. Set $H_t=K I_d-\Sigma_t$ and $h_-=K-\lambda_{\max}(\Sigma_0)>0$. Then $h_- I_d\preceq H_t\preceq K I_d$, since $\Sigma_t$ is a convex combination of $\Sigma_0$ and $(K/d)I_d$. Put $E_{a,b}=e_b e_a^T-e_a e_b^T$.
    The SDE equation
    \[
        \dd X_t=-2(d-1)\Lambda X_t\dd t
        +\sqrt{2\Lambda}\sum_{a<b}E_{a,b}X_t\dd\tilde{B}_t^{a,b}+\sqrt{2\Lambda}\,H_t^{1/2}\dd W_t,\quad
        \qquad X_0\sim f_0,
    \]
    with independent Brownian motions and initial value has mean zero and
    covariance $\Sigma_t$.  Moreover, its forward equation is the Maxwell Landau equation.

    Let $R_t\in\mathrm{SO}(d)$ solve
    \begin{equation}\label{eq:rotation-Strat}
        \mathrm dR_t=\sqrt{2\Lambda}
        \sum_{a<b}E_{a,b}R_t\circ\mathrm d\tilde{B}_t^{a,b},
        \qquad R_0=I.
    \end{equation}
    Skew-symmetry gives $R_t^TR_t=I$; continuity and $R_0=I$ give
    $\det R_t=1$. In It\^o form,
    \begin{equation}\label{eq:rotation-Ito}
        \mathrm dR_t=-(d-1)\Lambda R_t\dd t
        +\sqrt{2\Lambda}\sum_{a<b}E_{a,b}R_t\dd\tilde{B}_t^{a,b}.
    \end{equation}
    Variation of constants gives
    \begin{equation}\label{eq:X-representation}
        X_t=a_tR_tX_0+Z_t,
    \end{equation}
    where $a_t=e^{-(d-1)\Lambda t}$ and
    \begin{equation}\label{eq:Z}
        Z_t=\sqrt{2\Lambda}\,a_tR_t
        \int_0^t a_s^{-1}R_s^TH_s^{1/2}\dd W_s.
    \end{equation}
    Conditionally on the rotation path,
    $Z_t$ is Gaussian, independent of $X_0$, with covariance
    \[
        C_t=2\Lambda\int_0^t e^{-2(d-1)\Lambda(t-s)} R_tR_s^TH_sR_sR_t^T\dd s,
    \]
    with
    \[
        \frac{h_-}{d-1}(1-e^{-2(d-1)\Lambda t})I_d
        \le C_t\le
        \frac{K}{d-1}(1-e^{-2(d-1)\Lambda t})I_d.\]
    Consequently, with $\varphi_C$ the centered Gaussian density of covariance $C$
    \begin{equation*}
        \quad
        \varphi_C(z)
        =\frac{1}{(2\pi)^{d/2}\sqrt{\det C}}
        \exp\left(-\frac12 z^TC^{-1}z\right),
        \qquad z\in\R^d.
    \end{equation*}
    one may obtain that
    \begin{equation}\label{eq:step2-reference-representation}
        f(t,v)=\E\int f_0(y)\varphi_{C_t}(v-a_tR_ty)\dd y.
    \end{equation}
    Setting $\tilde{C}=a_t^{-2}R_t^TC_tR_t$ and using the formula of the Gaussian, one has
    \begin{gather}\label{eq:fformula}
        f(t, v)=\E a_t^{-d} f_0*\varphi_{\tilde{C}}(a_t^{-1}R_t^T v).
    \end{gather}
    Since $R_t$ is an orthogonal matrix and $a_t$ is bounded both below and above by positive numbers,
    one can find easily that $f\to f_0$ in $L^1$ and locally in $C^2$.
    This formula also proves positivity and positive-time smoothness.

    Using the explicit formula, one can directly derive that
    \begin{equation}\label{eq:first-density}
        \nabla_vf(t,v)
        =\E_R\int a_t^{-1}R_t\nabla f_0(y)
        \varphi_{C_t}(v-a_tR_ty)\dd y.
    \end{equation}
    and
    \begin{equation}\label{eq:second-density}
        \nabla_v^2f(t,v)
        =\E_R\int a_t^{-2}R_t\nabla^2f_0(y)R_t^T
        \varphi_{C_t}(v-a_tR_ty)\dd y.
    \end{equation}
    Let $\omega$ denote a rotation path and define
    \begin{equation}\label{eq:posterior}
        \nu_{t,v}(\dd\omega,\dd y)
        =\frac{f_0(y)\varphi_{C_t(\omega)}
        (v-a_tR_t(\omega)y)}{f(t,v)}
        \P(\dd\omega)\dd y.
    \end{equation}
    This is a conditional probability measure. Using this measure, one has
    \begin{equation}\label{eq:score-posterior}
        \nabla_v\log f(t,v)
        =\int a_t^{-1}R_t \nabla\log f_0(y)\dd\nu_{t,v},
    \end{equation}
    and
    \begin{equation}\label{eq:normalized-Hessian}
        \frac{\nabla_v^2f(t,v)}{f(t,v)}
        =\int a_t^{-2}R_t
        \bigl(\nabla^2\log f_0(y)+\nabla\log f_0(y)\otimes \nabla\log f_0(y)\bigr)
        R_t^T\dd\nu_{t,v}.
    \end{equation}

    We now make use of these explicit formulas to derive the desired bounds.
    The moment control follows directly from \Cref{eq:X-representation}.
    In fact, the first terms is easily controlled by $M_p(f_0)$ while the second term contributes
    \[
        \E |Z_t|^p \le C  M_2(f_0)^{p/2}\le C M_p(f_0).
    \]

    Next, we show that
    \[
        |\log f(t, v)|\le C_{T}\l v\r^{q_0}
    \]
    for $q_0=m_0+1\ge2$.
    First, the initial score bound and \Cref{eq:fformula} give the polynomial upper bounds
    \begin{gather}\label{eq:finfty}
        f_0(y)\le C\langle y\rangle^{d m_0},\qquad f(t,v)\le C_T\langle v\rangle^{d m_0}.
    \end{gather}
    Indeed, on the ball of radius $c\langle y\rangle^{-m_0}$ centered at $y$, the score bound implies $f_0(z)\ge c' f_0(y)$; integration over this ball proves the first inequality. The second follows from \Cref{eq:fformula} and Gaussian moments. The logarithmic derivative bounds also give polynomial bounds on the first two derivatives of $f_0$, justifying the differentiated Gaussian formulas and their local convergence at $t=0$.
    Moreover, integrating $\nabla\log f_0$ along a segment gives
    $f_0(v)\ge c e^{-C\langle v\rangle^{q_0}}$.
    For $0<t\le T$, the covariance bounds give $c_TtI_d\preceq C_t\preceq C_TtI_d$, uniformly in the rotation path.
    Restricting
    \eqref{eq:step2-reference-representation} to
    $|y-a_t^{-1}R_t^{\mathsf T}v|\le\sqrt t$ and using the uniform Gaussian
    bounds therefore yields
    \[
        f(t,v)\ge c_Te^{-C_T\langle v\rangle^{q_0}},\qquad 0<t\le T.
    \]
    Combining with \Cref{eq:finfty}, one has the desired control of $|\log f|$. The bounds also hold at $t=0$.

    Then, the estimates of $|\nabla\log f|^2$ and $|\nabla^2\log f|$ (equivalently, $|\nabla^2 f|/f$) will be similar to
    \cite[proof of Proposition~A.5]{li2026estimates}. In particular, both will be reduced to estimate of
    \[
        \int \l y\r^{p_0} d\nu_{t,v} =\E \int \l y\r^{p_0}
        \frac{f_0(y)\varphi_{C_t}(v-a_t R_t y)}{f(t, v)} dy.
    \]
    Here $p_0>0$. Since $R_t$ is orthogonal and $a_t$ is both bounded from above and below by positive constants, we first apply $\l y\r^{p_0}\le C(\l v\r^{p_0}+|v-a_t R_t y|^{p_0})$ and obtain
    \[
        \int \l y\r^{p_0} d\nu_{t,v}\le C \l v\r^{p_0}+CP, \quad P:=\E \int |v-a_t R_t y|^{p_0}
        \frac{f_0(y)\varphi_{C_t}(v-a_t R_t y)}{f(t, v)} dy.
    \]
    Write $z=v-a_tR_ty$. Taking $L=A\sqrt{t}\max\{1,\sqrt{|\log f(t,v)|+\log\langle v\rangle}\}$, one has
    \[
        P=P_1+P_2,
    \]
    where $P_1$ and $P_2$ are the integrals defining $P$ over $|z|\le L$ and $|z|>L$, respectively.
    One clearly has
    $P_1\le L^{p_0}$. For $P_2$, we split the Gaussian kernel with $C=C_t$ as
    \[
        \varphi_{C_t}(z)
        \le \frac{1}{(2\pi)^{d/2}\sqrt{\det C}}
        \exp\left(-\frac14 z^TC^{-1}z\right)\exp\left(-\frac18 z^TC^{-1}z\right)
        \exp\left(-\frac18\|C_t\|^{-1}|z|^2\right).
    \]
    We first observe that by choosing $A$ sufficiently large, one has
    \[
        \frac{\exp\left(-\frac18\|C_t\|^{-1}|v-a_tR_ty|^2\right)}{f(t,v)}
        \le C_T\langle v\rangle^{-d m_0},\qquad |z|>L.
    \]
    Indeed, $\|C_t\|\le C_Tt$, and the two terms inside the square root in $L$ absorb $1/f(t,v)$ and $\langle v\rangle^{d m_0}$, respectively. By \Cref{eq:finfty} and $\langle y\rangle\le C_T\langle v\rangle\langle z\rangle$, we also have $f_0(y)\le C_T\langle v\rangle^{d m_0}\langle z\rangle^{d m_0}$.
    Moreover,
    \[
        \sup_{z\in\R^d}|z|^{p_0}\langle z\rangle^{d m_0}
        \exp\left(-\frac18 z^TC_t^{-1}z\right)
        \le C_{T,p_0}t^{p_0/2}.
    \]
    Hence, one has
    \begin{multline*}
        P_2\le C_{T,p_0}t^{p_0/2}\\
        {}\times\E\int_{|v-a_tR_ty|>L}
        \frac{1}{(2\pi)^{d/2}\sqrt{\det C}}
        \exp\left(-\frac14(v-a_tR_ty)^TC^{-1}(v-a_tR_ty)\right)\,dy
        \le C_{T,p_0}.
    \end{multline*}
    The remaining Gaussian integral is bounded uniformly by the change of variables $z=v-a_tR_ty$.
    Combining the two, one thus has
    \[
        \int\l y\r^{p_0}d\nu_{t,v}
        \le C_{T,p_0}\bigl(\langle v\rangle^{p_0}+1+|\log f(t,v)|^{p_0/2}\bigr).
    \]
    With the bound of $|\log f|$, one has $\int\langle y\rangle^{p_0}d\nu_{t,v}\le C_{T,p_0}\langle v\rangle^{p_0q_0/2}$.
    Jensen's inequality in \Cref{eq:score-posterior} and the bound in \Cref{eq:normalized-Hessian}, with $p_0=2m_0$, yield
    \[
        |\nabla\log f|^2+\frac{|\nabla^2f|}{f}
        \le C_T\int\langle y\rangle^{2m_0}d\nu_{t,v}
        \le C_T\langle v\rangle^{m_0(m_0+1)}.
    \]
    Finally, $\nabla^2\log f=(\nabla^2f)/f-\nabla\log f\otimes\nabla\log f$ proves \Cref{eq:step2-reference-bound}.
\end{proof}

\section{Maxwell relative-coercivity estimate}\label{app:coercivity}

Here, we provide a reference proof for \Cref{lem:relativecoercivity}.

\begin{proof}
    Expanding the square in $D(\tilde{f}|f)$, one has
    \begin{equation}\label{eq:Dexpand}
        D(\tilde{f} | f)=\iint \tilde{f} \tilde{f}_* q(v)\cdot A(v-v_*)q(v)\dd v \dd v_*-\iint q(v)\cdot A(v-v_*) q(v_*)\tilde{f}\tilde{f}_*\dd v \dd v_*.
    \end{equation}
    Note
    \[
        A*\tilde{f}=A(v-v_c)+\Lambda((\tr \Sigma_s)I_d -\Sigma_s),
    \]
    where $\Sigma_s$ is introduced in \Cref{lem:covariance}. Note that $\tr(\Sigma_s)=\tr(\Sigma_0)=M_2-|v_c|^2
        =: K$.

    By \Cref{lem:covariance}, the covariance at any intermediate time has the form
    \[
        \Sigma(\tilde f(t))=\frac Kd I_d+\alpha(t)\left(\Sigma_0-\frac Kd I_d\right),\qquad 0\le\alpha(t)\le1.
    \]
    Here $\alpha(t)$ is the product of the factors $(1+e^{-8\Lambda d\tau})/2$ from completed steps and the analogous factor for the current step. Thus $\lambda_{\max}(\Sigma(\tilde f(t)))\le\lambda_{\max}(\Sigma_0)$, and
    \[
        A*\tilde f\succeq\Lambda\bigl(K-\lambda_{\max}(\Sigma_0)\bigr)I_d=\kappa_0I_d.
    \]
    Consequently,
    \[
        D(\tilde f\mid f)\ge\kappa_0\cI(\tilde f\mid f)-\mathcal C_f,\qquad
        \mathcal C_f=\iint q\cdot A q_*\tilde f\tilde f_*.
    \]

    To treat $\mathcal{C}_f$, we will make use of similar technique as in the proof of \Cref{lem:prooflem2}.
    The identity
    \begin{equation}\label{eq:gq}
        \tilde{f} q=\nabla(\tilde{f}-f)-(\tilde{f}-f)\nabla\log f
        =f\nabla(\tilde{f}/f-1)
    \end{equation}
    allows integration by parts in both velocities:
    \begin{equation}\label{eq:crossK}
        \mathcal C_f=\iint(\tilde{f}-f)(\tilde{f}_*-f_*)K^{\mathcal C}_f(v,v_*)\dd v\dd v_*,
    \end{equation}
    where with $b(z)=\nabla\cdot A(z)=-\Lambda(d-1) z$,
    \begin{equation}\label{eq:Kf}
        K^{\mathcal C}_f=-\nabla\cdot b(v-v_*)+b(\nabla \log f_*-\nabla\log f)
        +\nabla\log f(v)\cdot A(v-v_*)\nabla\log f(v_*).
    \end{equation}

    We recall the estimate of the reference solution in \Cref{lem:referencesol}, then
    \[
        |K^{\mathcal C}_f(v,w)|\le C_T\l v\r^{m+2}\l w\r^{m+2},
    \]
    with $m=m_0(m_0+1)$. Set $W(v):=\langle v\rangle^{m+2}$, so that $|K^{\mathcal C}_f(v,w)|\le C_TW(v)W(w)$. We recall the elementary square-root bound again
    \begin{equation}\label{eq:weightedHellinger}
        \left(\int W|g-f|\right)^2
        \le2H(g\mid f)\int W^2(g+f),\qquad W\ge0.
    \end{equation}
    where
    $\int(\sqrt g-\sqrt f)^2\le \cH(g\mid f)$ is used. Applying this inequality with $g=\tilde f$ to the kernel representation of $\mathcal C_f$, we obtain
    \[
        \begin{aligned}
            |\mathcal C_f|
             & \le C_T\left(\int W(v)|\tilde f(v)-f(v)|\dd v\right)^2                  \\
             & \le 2C_T\cH(\tilde f\mid f)\int\langle v\rangle^{2m+4}(\tilde f+f)\dd v \\
             & \le C\cH(\tilde f\mid f).
        \end{aligned}
    \]
    The last inequality uses the uniform $(2m+4)$-moment bounds for both $\tilde f$ and $f$, propagated under \Cref{ass:step2-initial}. Combining this estimate with the preceding lower bound for the dissipation yields
    \begin{equation}\label{eq:coercivity-appendix}
        D(\tilde{f}\mid f)\ge \kappa_0\cI(\tilde{f}\mid f)
        -C H(\tilde{f}\mid f).
    \end{equation}
    The weight on the right of \Cref{eq:weightedHellinger} requires moments of order $2m+4$, which are finite and uniformly propagated under \Cref{ass:step2-initial}. Spatial cutoffs justify both integrations by parts: their boundary terms tend to zero by these moment bounds and the weighted derivative bounds in \Cref{eq:pair-weighted-control}. This argument uses only first logarithmic derivatives of $f$.

\end{proof}

\section{Justification of the entropy identity}
\label[appendix]{app:entropy-justification}

We justify \Cref{eq:entropyderivative} without assuming differentiability of the untruncated entropy. Write $g=\tilde f$, $a=A*f$, $a_g=A*g$, and $b(v)=-\Lambda(d-1)(v-v_c)$. On each open time step,
\[
    \partial_t f=\nabla\cdot(a\nabla f-bf),\qquad
    \partial_t g=\nabla\cdot(a\nabla g-bg+V),\qquad
    V=E_{t-t_k}+(a_g-a)\nabla g.
\]
The two densities have the same mean and energy, so $a_g-a=-\Lambda(\Sigma_g-\Sigma_f)$ is independent of $v$ and uniformly bounded. The covariance formulas give $a\succeq\kappa_0I_d$, $|a|\le C_T\langle v\rangle^2$, and $|\nabla\cdot a|+|b|\le C_T\langle v\rangle$.

The propagation estimates and H\"older's inequality imply
\[
    \sup_{t\le T}\int\langle v\rangle^2\frac{|\nabla g|^2}{g}\dd v<\infty.
\]
Indeed, the moment required here is $2q/(q-2)<P_q$. Combining this with \Cref{lem:referencesol,eq:step2-uniform-Theta} yields
\begin{equation}\label{eq:entropy-integrability}
    \int_0^T\int\left[g\,q(v)\cdot a q(v)+\frac{|V|^2}{g}\right]\dd v\dd t<\infty,
    \qquad q(v)=\nabla\log(g/f).
\end{equation}
The contribution of the reference score in the first term requires only the moment of order $m+2\le p_*$. The pair derivative bounds give $\partial_t g,\nabla g,\nabla^2g,V,\nabla V\in L^1_{\mathrm{loc}}$ on each open step. Strong continuity of rotations on $L^1$, followed by expectation and marginalization, gives $g\in C([0,T];L^1)$; the endpoint values agree. The reference representation gives the same continuity for $f$.

Let $h=g/f$. Choose smooth functions $0\le\theta_M\uparrow1$ on $(0,\infty)$, supported in $[1/M,M]$, and define
\[
    \beta_M''(r)=\frac{\theta_M(r)}r,\qquad
    \beta_M(1)=\beta_M'(1)=0,\qquad \eta_M=f\beta_M(h).
\]
Then $\beta_M\uparrow r\log r-r+1$, while $\beta_M'$ and $\beta_M-r\beta_M'$ are bounded for fixed $M$. Thus $(f,g)\mapsto f\beta_M(g/f)$ is Lipschitz on the nonnegative quadrant, and $\eta_M$ is continuous in $L^1$ with $\eta_M(0)=0$ and $0\le\eta_M\le C_M(f+g)$. The local chain rule gives
\[
    \partial_t\eta_M
    =\nabla\cdot\bigl(a\nabla\eta_M-b\eta_M+\beta_M'(h)V\bigr)
    -\theta_M(h)\bigl[gq\cdot aq+q\cdot V\bigr].
\]
On compact sets at positive times, $f$ is smooth and bounded away from zero. Where $\beta_M''(h)\ne0$, $g/f$ is bounded above and below. The Sobolev chain rule is therefore justified by the local regularity above and the Fisher bound; all quadratic products are integrable by \Cref{eq:entropy-integrability}.

Test with a smooth spatial cutoff $\chi_R(v)=\chi(v/R)$, equal to one on $|v|\le R$ and zero on $|v|\ge2R$. Integrating the diffusion boundary flux once more by parts gives
\[
    \int\eta_M\bigl[\nabla\cdot(a\nabla\chi_R)+b\cdot\nabla\chi_R\bigr]
    -\int\beta_M'(h)V\cdot\nabla\chi_R.
\]
The bracket is bounded by $C_T\mathbf1_{\{R\le|v|\le2R\}}$. Hence its time integral vanishes as $R\to\infty$. The other term vanishes by Cauchy--Schwarz, \Cref{eq:entropy-integrability}, and $\int g=1$. At fixed $M$, $L^1$ continuity permits the lower time endpoint to tend to zero and cancels the contributions at mesh points. We obtain
\[
    \int\eta_M(t)=-\int_0^t\int\theta_M(h)\bigl[gq\cdot aq+q\cdot V\bigr]\dd v\dd u.
\]
Since $\int_0^T\int|q\cdot V|<\infty$, monotone convergence on the left and dominated convergence on the right give
\[
    H(t)=-\int_0^t\int\bigl[gq\cdot aq+q\cdot V\bigr]\dd v\dd u.
\]
In particular, $H$ is finite and absolutely continuous on $[0,T]$, with $H(0)=0$. Expanding $V$ and integrating by parts in the partner variable gives exactly $H'=-D+R_1+R_2$ almost everywhere. Those integrations by parts are justified by the same weighted derivative and moment bounds used in \Cref{app:coercivity}.

We also recover differentiability at every time inside each step. The orthogonal representation gives continuity of $F_s$ and its first derivatives in weighted $L^1$, hence continuity of $g$, $\nabla g$, and $E_s$ in $L^1$ on that step. Along any convergent time sequence, these fields therefore converge almost everywhere after passing to a subsequence. Positivity of $g$ and smoothness of $f$ at positive times give convergence of the integrand $gq\cdot aq+q\cdot V$. The higher-power derivative and flux bounds in \Cref{eq:step2-one-propagation,eq:local-three-norms}, together with the moment bound of order $p_*$, give uniform integrability on compact sets and uniformly vanishing spatial tails. Indeed, $q>6$, $r=4q/(q-2)>2$, and $p_*\ge2m+4$ leave strict room in H\"older's inequality for these quadratic products. Vitali's theorem then shows that the spatial integral is continuous on each open step. The preceding integral identity consequently gives $H\in C^1((t_k,t_{k+1}))$ for each $k=0,\ldots,K-1$, while $H$ remains continuous at the mesh points. This justifies applying Gr\"onwall's inequality separately on the time intervals.

\bibliographystyle{abbrv}
\bibliography{references}

\begin{thebibliography}{10}

\bibitem{AdamsEtAl2017}
M.~F. Adams, E.~Hirvijoki, M.~G. Knepley, J.~Brown, T.~Isaac, and R.~Mills.
\newblock {Landau} collision integral solver with adaptive mesh refinement on
  emerging architectures.
\newblock {\em SIAM Journal on Scientific Computing}, 39(6):C452--C465, 2017.

\bibitem{BailoCarrilloHu2024}
R.~Bailo, J.~A. Carrillo, and J.~Hu.
\newblock The collisional particle-in-cell method for the
  {Vlasov--Maxwell--Landau} equations.
\newblock {\em Journal of Plasma Physics}, 90(4):905900415, 2024.

\bibitem{BailoCarrilloMedagliaZanella2025}
R.~Bailo, J.~A. Carrillo, A.~Medaglia, and M.~Zanella.
\newblock Uncertainty quantification for the homogeneous
  {Landau--Fokker--Planck} equation via deterministic particle {Galerkin}
  methods.
\newblock {\em Multiscale Modeling \& Simulation}, 23(2):687--710, 2025.

\bibitem{Bird1994}
G.~A. Bird.
\newblock {\em Molecular Gas Dynamics and the Direct Simulation of Gas Flows}.
\newblock Clarendon Press, Oxford, 1994.

\bibitem{bobkov2001hypercontractivity}
S.~Bobkov, I.~Gentil, and M.~Ledoux.
\newblock Hypercontractivity of {Hamilton--Jacobi} equations.
\newblock {\em Journal de Math{\'e}matiques Pures et Appliqu{\'e}es},
  80(7):669--696, 2001.

\bibitem{BobylevNanbu2000}
A.~V. Bobylev and K.~Nanbu.
\newblock Theory of collision algorithms for gases and plasmas based on the
  {Boltzmann} equation and the {Landau--Fokker--Planck} equation.
\newblock {\em Physical Review E}, 61(4):4576--4586, 2000.

\bibitem{bolley2005weighted}
F.~Bolley and C.~Villani.
\newblock Weighted {Csisz{\'a}r-Kullback-Pinsker} inequalities and applications
  to transportation inequalities.
\newblock {\em Annales de la Facult{\'e} des sciences de Toulouse:
  Math{\'e}matiques}, 14(3):331--352, 2005.

\bibitem{BorghiPareschi2025}
G.~Borghi and L.~Pareschi.
\newblock {Wasserstein} convergence rates for stochastic particle approximation
  of {Boltzmann} models.
\newblock {\em SIAM Journal on Numerical Analysis}, 64(2):485--509, 2026.

\bibitem{Carrapatoso2016}
K.~Carrapatoso.
\newblock Propagation of chaos for the spatially homogeneous {Landau} equation
  for {Maxwellian} molecules.
\newblock {\em Kinetic and Related Models}, 9(1):1--49, 2016.

\bibitem{CarrilloDelgadinoDesvillettesWu2024}
J.~A. Carrillo, M.~G. Delgadino, L.~Desvillettes, and J.~S.-H. Wu.
\newblock The {Landau} equation as a gradient flow.
\newblock {\em Analysis \& PDE}, 17(4):1331--1375, 2024.

\bibitem{CarrilloDelgadinoWu2022Grazing}
J.~A. Carrillo, M.~G. Delgadino, and J.~Wu.
\newblock {Boltzmann} to {Landau} from the gradient flow perspective.
\newblock {\em Nonlinear Analysis}, 219:112824, 2022.

\bibitem{CFGJW}
J.~A. Carrillo, X.~Feng, S.~Guo, P.-E. Jabin, and Z.~Wang.
\newblock Relative entropy method for particle approximation of the {Landau}
  equation for {Maxwellian} molecules.
\newblock {\em Journal de Math{\'e}matiques Pures et Appliqu{\'e}es},
  206:103838, 2026.

\bibitem{carrillo2025fisherinformationdecaykac}
J.~A. Carrillo and S.~Guo.
\newblock From {F}isher information decay for the {K}ac model to the
  {L}andau-{C}oulomb hierarchy.
\newblock {\em Preprint arXiv:2502.18606}, 2025.

\bibitem{CarrilloGuoJabin2024}
J.~A. Carrillo, S.~Guo, and P.-E. Jabin.
\newblock Mean-field derivation of {Landau}-like equations.
\newblock {\em Applied Mathematics Letters}, 158:109195, 2024.

\bibitem{CarrilloHuVanFleet2024}
J.~A. Carrillo, J.~Hu, and S.~Q. Van~Fleet.
\newblock A particle method for the multispecies {Landau} equation.
\newblock {\em Acta Applicandae Mathematicae}, 194:6, 2024.

\bibitem{CarrilloHuWangWu2020}
J.~A. Carrillo, J.~Hu, L.~Wang, and J.~Wu.
\newblock A particle method for the homogeneous {Landau} equation.
\newblock {\em Journal of Computational Physics: X}, 7:100066, 2020.

\bibitem{CarrilloJinTang2022}
J.~A. Carrillo, S.~Jin, and Y.~Tang.
\newblock Random batch particle methods for the homogeneous {Landau} equation.
\newblock {\em Communications in Computational Physics}, 31(4):997--1019, 2022.

\bibitem{DMN}
F.~Delarue, S.~Menozzi, and E.~Nualart.
\newblock The {Landau} equation for {Maxwellian} molecules and the brownian
  motion on $\mathrm{SO}_{N}(\mathbb{R})$.
\newblock {\em Electronic Journal of Probability}, 20(92):1--39, 2015.

\bibitem{DimarcoLiPareschiYan2015}
G.~Dimarco, Q.~Li, L.~Pareschi, and B.~Yan.
\newblock Numerical methods for plasma physics in collisional regimes.
\newblock {\em Journal of Plasma Physics}, 81(1):305810106, 2015.

\bibitem{du2024collision}
K.~Du and L.~Li.
\newblock A collision-oriented interacting particle system for {Landau}-type
  equations and the molecular chaos.
\newblock {\em arXiv preprint arXiv:2408.16252}, 2024.

\bibitem{du2025structure}
K.~Du, L.~Li, Y.~Xie, and Y.~Yu.
\newblock A structure-preserving collisional particle method for the {Landau}
  kinetic equation.
\newblock {\em Journal of Computational Physics}, 543:114387, 2025.

\bibitem{FengWang2026Kac}
X.~Feng and Z.~Wang.
\newblock {Kac}'s program for the {Landau} equation.
\newblock {\em Archive for Rational Mechanics and Analysis}, 250(4):72, 2026.

\bibitem{Fournier2009}
N.~Fournier.
\newblock Particle approximation of some {Landau} equations.
\newblock {\em Kinetic and Related Models}, 2(3):451--464, 2009.

\bibitem{FournierGuillin2015}
N.~Fournier and A.~Guillin.
\newblock On the rate of convergence in {Wasserstein} distance of the empirical
  measure.
\newblock {\em Probability Theory and Related Fields}, 162(3--4):707--738,
  2015.

\bibitem{FGHard}
N.~Fournier and A.~Guillin.
\newblock From a {Kac}-like particle system to the {Landau} equation for hard
  potentials and {Maxwell} molecules.
\newblock {\em Annales scientifiques de l'{\'E}cole Normale Sup{\'e}rieure},
  50:157--199, 2017.

\bibitem{GuillenSilvestre2025}
N.~Guill{\'e}n and L.~Silvestre.
\newblock The {Landau} equation does not blow up.
\newblock {\em Acta Mathematica}, 234:315--375, 2025.
\newblock Pair estimates are referenced to arXiv:2311.09420v2.

\bibitem{HirvijokiAdams2017}
E.~Hirvijoki and M.~F. Adams.
\newblock Conservative discretization of the {Landau} collision integral.
\newblock {\em Physics of Plasmas}, 24(3):032121, 2017.

\bibitem{JinLi2022MeanFieldRBM}
S.~Jin and L.~Li.
\newblock On the mean field limit of the {Random Batch Method} for interacting
  particle systems.
\newblock {\em Science China Mathematics}, 65(1):169--202, 2022.

\bibitem{JinLiLiu2020}
S.~Jin, L.~Li, and J.-G. Liu.
\newblock Random batch methods ({RBM}) for interacting particle systems.
\newblock {\em Journal of Computational Physics}, 400:108877, 2020.

\bibitem{JinLiLiu2021Disparate}
S.~Jin, L.~Li, and J.-G. Liu.
\newblock Convergence of the random batch method for interacting particles with
  disparate species and weights.
\newblock {\em SIAM Journal on Numerical Analysis}, 59(2):746--768, 2021.

\bibitem{Kallenberg2021}
O.~Kallenberg.
\newblock {\em Foundations of Modern Probability}, volume~99 of {\em
  Probability Theory and Stochastic Modelling}.
\newblock Springer, Cham, third edition, 2021.

\bibitem{Landau1936}
L.~D. Landau.
\newblock Die kinetische gleichung f{\"u}r den fall {Coulombscher}
  wechselwirkung.
\newblock {\em Physikalische Zeitschrift der Sowjetunion}, 10:154--164, 1936.

\bibitem{li2026estimates}
L.~Li, M.~Wang, and Y.~Wang.
\newblock Estimates of the numerical density for stochastic differential
  equations with multiplicative noise.
\newblock {\em Science China Mathematics}, 69(8):2127--2160, 2026.

\bibitem{MedagliaPareschiZanella2024}
A.~Medaglia, L.~Pareschi, and M.~Zanella.
\newblock Particle simulation methods for the {Landau--Fokker--Planck} equation
  with uncertain data.
\newblock {\em Journal of Computational Physics}, 503:112845, 2024.

\bibitem{Nanbu1980}
K.~Nanbu.
\newblock Direct simulation scheme derived from the {Boltzmann} equation. {I}.
  monocomponent gases.
\newblock {\em Journal of the Physical Society of Japan}, 49(5):2042--2049,
  1980.

\bibitem{PareschiRussoToscani2000}
L.~Pareschi, G.~Russo, and G.~Toscani.
\newblock Fast spectral methods for the {Fokker--Planck--Landau} collision
  operator.
\newblock {\em Journal of Computational Physics}, 165(1):216--236, 2000.

\bibitem{tabary2026weak}
C.~Tabary.
\newblock Weak-strong uniqueness for the {Landau} equation by a relative
  entropy method.
\newblock {\em SIAM Journal on Mathematical Analysis}, 58(5):4601--4640, 2026.

\bibitem{talagrand1996transportation}
M.~Talagrand.
\newblock Transportation cost for {Gaussian} and other product measures.
\newblock {\em Geometric \& Functional Analysis GAFA}, 6(3):587--600, 1996.

\bibitem{Trevisan2016}
D.~Trevisan.
\newblock Well-posedness of multidimensional diffusion processes with weakly
  differentiable coefficients.
\newblock {\em Electronic Journal of Probability}, 21(22):1--41, 2016.

\bibitem{vanErvenHarremoes2014}
T.~van Erven and P.~Harremo{\"e}s.
\newblock {R{\'e}nyi} divergence and {Kullback--Leibler} divergence.
\newblock {\em IEEE Transactions on Information Theory}, 60(7):3797--3820,
  2014.

\bibitem{Villani2002}
C.~Villani.
\newblock A review of mathematical topics in collisional kinetic theory.
\newblock In S.~Friedlander and D.~Serre, editors, {\em Handbook of
  Mathematical Fluid Dynamics}, volume~1, pages 71--305. North-Holland,
  Amsterdam, 2002.

\bibitem{Villani2009}
C.~Villani.
\newblock {\em Optimal Transport: Old and New}.
\newblock Springer, Berlin, 2009.

\end{thebibliography}

\end{document}